\documentclass{cmslatex}
\usepackage[paperwidth=7in, paperheight=10in, margin=.875in]{geometry}
\usepackage[backref,colorlinks,linkcolor=red,anchorcolor=green,citecolor=blue]{hyperref}
\usepackage{amsfonts,amssymb}
\usepackage{amsmath}
\usepackage{graphicx}
\usepackage{cite}
\usepackage{enumerate}

\graphicspath{{fig/}}

\renewcommand{\theequation}{\arabic{section}.\arabic{equation}}

	\DeclareMathOperator*{\argmin}{argmin}

	\newcommand{\br}{\mathbb{R}}
	\newcommand{\rd}{\mathrm{d}}
	\newcommand{\mcd}{\mathcal{D}}

	\newcommand{\boeta}{\boldsymbol{\eta}}
	\newcommand{\bn}{\boldsymbol{n}}
	\newcommand{\bx}{\boldsymbol{x}}
	\newcommand{\bxi}{\boldsymbol{\xi}}

	\newcommand{\abs}[1]{\left\vert#1\right\vert}
	\newcommand{\brac}[1]{\left(#1\right)}
	\newcommand{\cbra}[1]{\left\{#1\right\}}

	\newcommand{\WFR}{\textrm{WFR}_{\gamma}}	
	
	\def\half{\frac 1 2}

   \allowdisplaybreaks

\begin{document}
	\title{The Wasserstein-Fisher-Rao metric for waveform based earthquake location\thanks{Received date, and accepted date (The correct dates will be entered by the editor).} This work was supported by the National Natural Science Foundation of China (Grant Nos. 11871297, 41390452, 91730306) and the National Key R\&D Program on Monitoring, Early Warning and Prevention of Major Natural Disaster (Grant No. 2017YFC1500301).}
    
	%For each author, make a block with the following macros:

	\author{Datong Zhou\thanks{Department of Mathematical Sciences, Tsinghua University, Beijing, China 100084 (zdt14@mails.tsinghua.edu.cn).}
		\and Jing Chen\thanks{Department of Mathematical Sciences, Tsinghua University, Beijing, China 100084 (jing-che16@mails.tsinghua.edu.cn).}
		\and Hao Wu\thanks{Corresponding author. Department of Mathematical Sciences, Tsinghua University, Beijing, China 100084 (hwu@tsinghua.edu.cn).}
		\and Dinghui Yang\thanks{Department of Mathematical Sciences, Tsinghua University, Beijing, China 100084 (dhyang@math.tsinghua.edu.cn).}
		\and Lingyun Qiu\thanks{Yau Mathematical Sciences Center, Tsinghua University, Beijing China 100084 (Qiu.Lingyun@gmail.com).}}

	\pagestyle{myheadings} \markboth{The Wasserstein-Fisher-Rao metric for waveform based earthquake location}{D.T. Zhou, J. Chen, H. Wu, D.H. Yang and L.Y. Qiu}
	
    \maketitle

	\begin{abstract}
		In our previous work [Chen el al., J. Comput. Phys., 373(2018)], the quadratic Wasserstein metric is successfully applied to the earthquake location problem. The actual earthquake hypocenter can be accurately recovered  starting from initial values very far from the true ones. However, the seismic wave signals need to be normalized since the quadratic Wasserstein metric requires mass conservation. This brings a critical difficulty. Since the amplitude of a seismogram at a receiver is a good representation of the distance between the source and the receiver, simply normalizing the signals will cause the objective function in optimization process to be insensitive to the distance between the source and the receiver. When the data is contaminated with strong noise, the minimum point of the objective function will deviate and lead to a low accurate location result.
		
		To overcome the difficulty mentioned above, we apply the Wasserstein-Fisher-Rao (WFR) metric [Chizat et al., Found. Comput. Math., 18(2018)] to the earthquake location problem. The WFR metric is one of the newly developed metric in the unbalanced Optimal Transport theory. It does not require the normalization of the seismic signals. Thus, the amplitude of seismograms can be considered as a new constraint, which can substantially improve the sensitivity of the objective function to the distance between the source and the receiver. As a result, we can expect more accurate location results from the WFR metric based method compare to those based on quadratic Wasserstein metric under high-intensity noise. The numerical examples also demonstrate this.
	\end{abstract}

	\begin{keywords}
		the Wasserstein-Fisher-Rao metric; the quadratic Wasserstein metric; Inverse theory; Waveform inversion; Earthquake location
	\end{keywords}

	\begin{AMS} 49N45; 65K10; 86-08; 86A15 \end{AMS}

% ********************************************************************************
% ************* Introduction *****************************************************
% *******************************************************************************
\section{Introduction}\label{sec:intro}
The waveform based earthquake location models \cite{ChChWuYa:18, KiLiTr:11, LiPoKoTr:04, ToZhYaYaChLi:14, WuChHuYa:16, WuChJiToYa:17} have attracted wide attention from both academia and industry due to their highly accurate  results in deriving the hypocenter $\bxi_T$ and origin time $\tau_T$ of an earthquake, which are very important for the accurate inversion of underground velocity structure \cite{WaEl:00}. In addition, the precise locations of a large number of small-scale earthquakes in a particular area allow the seismogenic structure in the area to be obtained. Thus, the early warning predictions of an earthquake may be made according to the investigation of the structure \cite{SaLoZo:08}. Moreover, by determining the location of micro-earthquakes, we can monitor the dynamic disaster in mines and hydro-fracture \cite{EaFo:11, LeSt:81}.

From the mathematical point of view, the waveform based earthquake location problem  can be written as a nonlinear optimization model with PDE constraints \cite{ChChWuYa:18, WuChHuYa:16, WuChJiToYa:17}
\begin{equation} \label{eqn:inv_prob}
	(\bxi_T,\tau_T)=\argmin_{\bxi,\tau}\sum_r\chi_r(\bxi,\tau).
\end{equation}
Here $r$ denotes the index of receivers, and the corresponding misfit function $\chi_r(\bxi,\tau)$ is defined as
\begin{equation} \label{eqn:mist_fun}
	\chi_r(\bxi,\tau)=\mcd^2 (d_r(t),s(\boeta_r,t)).
\end{equation}
In the above equation, $d_r(t)$ and $s(\bx,t)$ are the real and synthetic earthquake signals respectively, and $\mcd$ measures the distance between them, which will be specified later. The two wavefields, 
\begin{equation} \label{eqn:rel_syn}
	d_r(t)=u(\boeta_r,t;\bxi_T,\tau_T), \quad s(\bx,t)=u(\bx,t;\bxi,\tau),
\end{equation}
can be regarded as the solutions to the following acoustic wave equation with initial-boundary conditions
\begin{align}
	& \frac{\partial^2 u(\bx,t;\bxi,\tau)}{\partial t^2}=\nabla\cdot\brac{c^2(\bx)\nabla u(\bx,t;\bxi,\tau)}
		+R(t-\tau)\delta(\bx-\bxi), \quad \bx,\bxi\in\Omega, \label{eqn:wave} \\
	& u(\bx,0;\bxi,\tau)=\partial_t u(\bx,0;\bxi,\tau)=0, \quad \bx\in\Omega, \label{ic:wave} \\
	& \bn\cdot \brac{c^2(\bx)\nabla u(\bx,t;\bxi,\tau)}=0, \quad \bx\in\partial \Omega. \label{bc:wave}
\end{align}
Here $c(\bx)$ is the acoustic wave speed and $\boeta_r$ denotes the location of the $r-$th receiver. The point source $\delta(\bx-\bxi)$ is used to model the seismic rupture inside the computational domain $\Omega\subset \br^d$. This assumption is reasonable since the scale of seismic rupture is much smaller than that of the seismic wave \cite{Ma:15}. The Ricker wavelet is used as the source wavelet
\begin{equation} \label{eqn:ricker}
	R(t)=A(1-2\pi^2f_0^2t^2)e^{-\pi^2f_0^2t},
\end{equation}
in which $f_0$ is the dominant frequency and $A$ is the indication parameter of amplitude. On the boundary $\partial \Omega$, the reflection condition \eqref{bc:wave} is used for simplification and $\bn$ is the outward unit normal vector to the domain $\Omega$. We can also easily consider the other boundary conditions, e.g., the perfectly matched layer absorbing boundary condition \cite{KoTr:03}.

Compared to the ray-based earthquake location methods \cite{Ge:03, Ge:03b, WaEl:00, WuYa:13}, the most significant advantage of the waveform based methods is that the location accuracy is much higher and the resolution can reach the wavelength scale \cite{RaPoFi:10}. However, some challenges are to be overcome. First, the waveform based earthquake location methods are computationally expensive since we need to solve the wave equation many times to obtain the synthetic signals and construct the constraints for inversion. Second, the waveform inversion with $\ell^2$ norm is suffering from the well-known cycle-skipping phenomenon, which is associated with many local minima of the objective function in the optimization. Third, the location results may be inaccurate when the real seismic signals are highly noisy. 

In recent years, the Wasserstein metric in the Optimal Transport theory was successfully applied to solve the inverse problems in seismology \cite{EnFr:14, EnFrYa:16, QiJaAlYaEn:17,MeBrMeOuVi:16, MeBrMeOuVi:16b, YaEnSuFr:18, YaEn:18}. The newly defined model significantly mitigates the difficulty of many local minima of the objective function. In addition, the Wasserstein metric is less sensitivity to the data noise. Thus, we can expect reasonably accurate inversion results in the cases when the data is contaminated with high-intensity noise. However, the quadratic Wasserstein metric also has limitations: it requires the signals to have the same integral value \cite{Vi:03, Vi:08}. In previous studies \cite{ChChWuYa:18, EnFr:14, EnFrYa:16, QiJaAlYaEn:17}, different normalization procedures are applied to the signals to address the issue. Such an operation ignores the important amplitude information which is a good representation of the distance between the source and the receiver. This may cause nearly flat objective function along certain direction since the distance between the source and receiver affect the arrival time of a signal at a receiver. Thus, the minimum point of the objective function may deviate a lot even under the small magnitude of data noise, which leads to low accurate location results, see Examples \ref{exam:two_layler_convexity}-\ref{exam:two_layler_noise_convexity_and_scatter} for illustration.

\vskip2mm\begin{remark}
	The Kantorovich-Rubinstein (KR) norm \cite{MeBrMeOuVi:16, MeBrMeOuVi:16b} does not require the signals to have the same integral. However, the convexity of the objective functions defined with the KR norm may not be good enough for the earthquake location problems, see Figures 2-3 in \cite{ChChWuYa:18} for illustration.
\end{remark}\vskip2mm

The Wasserstein-Fish-Rao (WFR) metric is a newly developed transport metric \cite{ChPeScVi:18, ChPeScVi:18b}, it allows to compute the distance between arbitrary positive signals, see also related topics in \cite{KoMoVo:16, LiMiSa:18, PiRo:14, PiRo:16}. This metric is an interpolation between the quadratic Wasserstein metric and the Fisher-Rao metric. From the fluid dynamics point of view \cite{BeBr:00}, the new metric introduces a source term in the continuity equation. Since it allows to measure the distance between two signals with different total integral, the normalization to the seismic signals is no longer required. Thus, the important amplitude information is retained, and the new objective function defined by this metric has better convexity near the exact earthquake hypocenter. Therefore, we can expect more accurate location results even with strong noise in the data, which is the main goal of our paper.

This paper is organized as follows. In Section \ref{sec:WFR}, the formulation, basic properties and the computational methods of the Wasserstein-Fisher-Rao metric are reviewed. We then apply this metric to the earthquake location problems in Section \ref{sec:ael}. In Section \ref{sec:numer}, the numerical experiments are provided to demonstrate the effectiveness and efficiency of the new method. Finally, we make some conclusive remarks in Section \ref{sec:cons}.

% ********************************************************************************
% ************ Wasserstein-Fisher-Rao metric *******************************
% ********************************************************************************
\section{The Wasserstein-Fisher-Rao metric}\label{sec:WFR}
For two non-negative distribution functions $\mu, \nu$ on a domain $\Omega \in \mathbb{R}^n$ and interpolating parameter $0<\gamma<\infty$, the WFR metric is defined by solving a minimizing problem:
\begin{equation} \label{eqn:WFR}
	\WFR(\mu,\nu) = \brac{ \inf_{\rho,v,\alpha} \int_0^1 \int_{\Omega}  \brac{\half \abs{v(t,x)}^2 + \frac{\gamma^2}{2} \alpha(t,x)^2 }\rho(t,x) \rd x \rd t}^{\half} ,
\end{equation}
under the constraint that the triplet $(\rho, v, \alpha)$ satisfy the following continuity equation:
\begin{equation} \label{eqn:WFR_con}
	\left\{\begin{array}{l}
		\partial_t \rho + \nabla \cdot(\rho v) = \rho \alpha , \\
		\rho(0,\cdot)=\mu, \quad \rho(1,\cdot) =\nu.
	\end{array}\right.
\end{equation}
Here $\rho$ is arbitrary time-dependent density, $v$ is arbitrary velocity field that stands for the movement of mass, and $\alpha$ is arbitrary scalar field associated with the creation and destruction of mass. This metric was defined and studied simultaneously and independently in \cite{ChPeScVi:18, KoMoVo:16, LiMiSa:16}, with quite different approaches. For $\abs{\mu}=\abs{\nu}$ and $\gamma \to \infty$, the source term $\rho \alpha$ in the continuity equation is depleted (i.e. $\alpha \equiv 0$) and the metric degenerate to the dynamical formulation of Wasserstein metric of Benamou and Brenier in \cite{BeBr:00}.
Developing the numerical method of the classical optimal transport problems is still a very active topic recently, and several numerical methods  are suggested \cite{BeBr:00, BeCaCuNePe:15, EnFrYa:16, Fr:12, Wa:04}. Some of these methods can be extended to the generalized problems introduced in \eqref{eqn:WFR}-\eqref{eqn:WFR_con}.

One famous numerical method under the continuous framework is the one proposed by Benamou and Brenier \cite{BeBr:00} that solves the convex variational problem. While the original paper only suggested the algorithm for classical Wasserstein metrics (the above minimizing problem restricted on $\gamma=\infty,\;\alpha\equiv0$), this framework can be adapted to dynamic terms modified in various ways, especially $\gamma<\infty$ that appears in WFR model. However, for an optimal transport problem in d-dimensional space, this framework requires solving the elliptic equation in a $(d+1)$-dimension space iteratively. This leads to a large computational cost, which makes the method not suitable for large-scale inversion problems.

Recently, a simple and efficient approach is proposed based on the linear programming interpretation of optimal transport problem. This approach makes use of the discretization of general non-negative distributions $\mu$ and $\nu$ by approximating them with two finite atomic distributions $\sum_{i = 1}^N a_i \delta_{x_i}$ and $\sum_{j = 1}^M b_j \delta_{y_j}$. Entropy regularization is then introduced to make the obtained discrete linear programming problem compatible with the celebrated Sinkhorn algorithm \cite{BeCaCuNePe:15, ChPeScVi:18, ChPeScVi:18b, Si:64}. As there exists a linear programming interpretation of WFR model, naturally this approach can be applied to the problem \eqref{eqn:WFR}-\eqref{eqn:WFR_con}. In the following, we will review the basic ideas in a heuristic way. For rigorous analysis, the readers are referred to \cite{LiMiSa:18}. We start with a lemma about the simplest case of transporting one Dirac function to another, the detailed proof can be found in Section 4.2 of \cite{ChPeScVi:18}. \smallskip

\begin{lemma}
	Let $\mu=h_0 \delta_{x_0}$ and $\nu = h_1 \delta_{x_1}$ be the initial and final distributions, where $x_0,\; x_1 \in \Omega$ are the locations and $h_0, h_1 \geq 0$ are the mass of Diracs that may be zero. Then in the Wasserstein-Fisher-Rao space, depending on the distance $\abs{x_1 - x_0}$, the geodesic between the two Diracs will behave in 3 distinct ways:
	\begin{enumerate}
		\item Traveling Dirac. If $\abs{x_1 - x_0} < \pi \gamma$, then the traveling Dirac
			\begin{equation*}
				\rho(t) = h(t) \delta_{x(t)},
			\end{equation*}
			implicitly defined by
			\begin{equation*}
				\left\{ \begin{aligned}
					& h(t) = At^2 - 2Bt +h_0, \\
					& h(t)x'(t) = \omega_0,
				\end{aligned} \right.
			\end{equation*}
			with
			\begin{align*}
				\omega_0 = 2 \gamma \tau \sqrt{\frac{h_0 h_1}{1 + \tau^2}}, \quad \tau = \tan \left( \frac{x_1 - x_0}{2\gamma} \right), \\
				A = h_1 + h_0 - 2 \sqrt{\frac{h_0h_1}{1 + \tau^2}}, \quad B = h_0 - \sqrt{\frac{h_0h_1}{1 + \tau^2}},
			\end{align*}
			is the unique geodesic.

		\item Cut Locus. If $\abs{x_1 - x_0} = \pi \gamma$, there are infinitely many geodesics. Two of the particular examples are traveling Dirac and Fisher-Rao geodesic.

		\item Fisher-Rao Geodesic (No Transport). If $\abs{x_1 - x_0} > \pi \gamma$, then the Fisher-Rao geodesic
			\begin{equation*}
				\rho(t) = t^2 h_1 \delta_{x_1} + (1 - t)^2 h_0 \delta_{x_0},
			\end{equation*}
			is the unique geodesic.
	\end{enumerate}
	As a consequence, the distance between two Diracs are
	\begin{equation} \label{eqn:WFR_Dirac}
		\WFR(h_0\delta_{x_0}, h_1\delta_{x_1}) = \sqrt{2} \gamma \left[h_0 + h_1 - 2\sqrt{h_0 h_1} \cos_+\brac{\frac{x_1 - x_0}{2\gamma}}\right]^{\half},
	\end{equation}
	where
	\begin{equation*}
		\cos_+\left( x \right) = \left\{ \begin{aligned}
			\cos(x) \quad x \in \left[-\frac{\pi}{2}, \frac{\pi}{2}\right], \\
			0 \quad x \notin \left[-\frac{\pi}{2}, \frac{\pi}{2}\right].
		\end{aligned} \right.
	\end{equation*}
\end{lemma}

\vskip2mm\begin{remark} \label{rem:gamma}
	The parameter $\gamma$ controls the interpolation between two terms of the minimization problem \eqref{eqn:WFR}. As the lemma has demonstrated, there will be coupling between two Diracs only if the distance between them is closer than $\pi \gamma$, yet changing $\gamma$ is equivalent to changing the scale of the problem. Intuitively, the larger $\gamma$ is, the lower amount of mass will be created or removed, and the behavior of WFR metric will be more like the classical Wasserstein metric. On the other hand, as $\gamma$ vanishes,  the behavior of WFR metric will approach Fisher-Rao metric. 
\end{remark}\vskip2mm

According to the above results, we can see the WFR metric \eqref{eqn:WFR_Dirac} is convex with respect to the mass of Diracs $h_0$ and $h_1$. Next, we consider a more general transport situation that the two atomic distributions are the summation of many Diracs
\begin{equation} \label{eqn:atomic_distribution}
	\mu = \sum_{i = 1}^N \mu_i \delta_{x_i}, \quad
	\nu = \sum_{j = 1}^M \nu_j \delta_{y_j}, \quad \mu_i\ge0,\;\nu_j\ge0.
\end{equation}
Naturally, we can split the mass $\mu_i,\;\nu_j$ into different pieces $\alpha_{ij}\ge0,\;\beta_{ji}\ge0$ as
\begin{align} 
	&\sum_{j = 1}^M \alpha_{ij} = \mu_i, \quad i=1,2,\cdots,N, \label{eqn:split_constraint1} \\
	& \sum_{i = 1}^N \beta_{ji} = \nu_j, \quad j=1,2,\cdots,M. \label{eqn:split_constraint2}
\end{align}
Obviously, we have
\begin{equation*}
	\WFR^2\brac{\sum_{i = 1}^N \mu_i \delta_{x_i},\sum_{j = 1}^M \nu_j \delta_{y_j}}
		\le \inf_{\alpha_{ij},\beta_{ji}}\sum_{i,j}\WFR^2\brac{\alpha_{ij}\delta_{x_i},\beta_{ji}\delta_{y_j}}.
\end{equation*}
On the other hand, since the functional is convex and homogeneous, we can also assert that the local (infinitesimal) behavior of the optimal solution must comply with the optimal solution in the local (infinitesimal) problem \cite{LiMiSa:18}. As a consequence,
\begin{equation*}
	\WFR^2\brac{\sum_{i = 1}^N \mu_i \delta_{x_i},\sum_{j = 1}^M \nu_j \delta_{y_j}}
		\ge \inf_{\alpha_{ij},\beta_{ji}}\sum_{i,j}\WFR^2\brac{\alpha_{ij}\delta_{x_i},\beta_{ji}\delta_{y_j}}.
\end{equation*}
Thus, the searching for optimal transport cost of the two atomic distributions $\mu$ and $\nu$ in \eqref{eqn:atomic_distribution} is equivalent to minimizing the total cost of 
\begin{equation} \label{eqn:WFR_atomic}
	\sum_{i,j}\WFR^2\brac{\alpha_{ij}\delta_{x_i},\beta_{ji}\delta_{y_j}}=
		2\gamma^2\sum_{i,j}\brac{\alpha_{ij}+\beta_{ji}-2\sqrt{\alpha_{ij}\beta_{ji}}\cos_+\brac{\frac{x_i-y_j}{2\gamma}}},
\end{equation}
under constraints \eqref{eqn:split_constraint1}-\eqref{eqn:split_constraint2}. Since the above function is convex with respect to the freedoms $\alpha_{ij}$ and $\beta_{ji}$ and bounded from below, the existence of the minimum is guaranteed. Combine \eqref{eqn:WFR_atomic} and \eqref{eqn:split_constraint1}-\eqref{eqn:split_constraint2} into one formula, we obtain
{\small\begin{multline}
	\WFR^2\brac{\sum_{i = 1}^N \mu_i \delta_{x_i},\sum_{j = 1}^M \nu_j \delta_{y_j}} \\
		=2\gamma^2\inf_{\alpha_{ij},\beta_{ji}}\cbra{\WFR^2\brac{\alpha_{ij}\delta_{x_i},\beta_{ji}\delta_{y_j}} \;:\;
			\sum_{j = 1}^M \alpha_{ij} = \mu_i,\;\sum_{i = 1}^N \beta_{ji} = \nu_j,\;\alpha_{ij}\ge0,\;\beta_{ji}\ge0}, \\
		=2\gamma^2\inf_{\alpha_{ij},\beta_{ji}\ge0}\Bigg\{\sum_{i,j}\brac{\alpha_{ij}+\beta_{ji}-2\sqrt{\alpha_{ij}\beta_{ji}}\cos_+\brac{\frac{x_i-y_j}{2\gamma}}} \\
			+\sup_{\phi_i}\cbra{\sum_i\phi_i\brac{\mu_i-\sum_j\alpha_{ij}}}
			+\sup_{\psi_j}\cbra{\sum_j\psi_j\brac{\nu_j-\sum_i\beta_{ji}}}\Bigg\} \\
		=2\gamma^2\inf_{\alpha_{ij},\beta_{ji}\ge0}\sup_{\phi_i,\psi_j}\Bigg\{\sum_i\phi_i\mu_i+\sum_j\psi_j\nu_j \\
			+\sum_{i,j}\brac{\brac{1-\phi_i}\alpha_{ij}+\brac{1-\psi_j}\beta_{ji}-2\sqrt{\alpha_{ij}\beta_{ji}}\cos_+\brac{\frac{x_i-y_j}{2\gamma}}}\Bigg\} \\
		=2\gamma^2\sup_{\phi_i,\psi_j}\inf_{\alpha_{ij},\beta_{ji}\ge0}\Bigg\{\sum_i\phi_i\mu_i+\sum_j\psi_j\nu_j \\
			+\sum_{i,j}\brac{\brac{1-\phi_i}\alpha_{ij}+\brac{1-\psi_j}\beta_{ji}-2\sqrt{\alpha_{ij}\beta_{ji}}\cos_+\brac{\frac{x_i-y_j}{2\gamma}}}\Bigg\} \\
		=2\gamma^2\sup_{\phi_i,\psi_j}\cbra{\sum_i\phi_i\mu_i+\sum_j\psi_j\nu_j \;:\;
			\phi_i,\psi_j\le 1, \; (1-\phi_i)(1-\psi_j)\ge \cos_+^2\brac{\frac{x_i-y_j}{2\gamma}}} \label{eqn:WFR_atomic_con}
\end{multline}}
The last second equality is an interchange of infimum and supremum that requires strong duality. Indeed, since $\WFR^2\brac{\alpha_{ij}\delta_{x_i},\beta_{ji}\delta_{y_j}}$ is convex with respect to $\alpha_{ij}, \beta_{ji}$ and we have assumed the space $X$ and $Y$ to be discrete and finite, strong duality and the existence of a minimizer is guaranteed by the Fenchel-Rockafellar Theorem \cite{Ro:67}. The constraints in the last expression guarantee that the lower bound of the functions
\begin{equation*}
	\Upsilon(\alpha_{ij},\beta_{ji})=\brac{1-\phi_i}\alpha_{ij}+\brac{1-\psi_j}\beta_{ji}-2\sqrt{\alpha_{ij}\beta_{ji}}\cos_+\brac{\frac{x_i-y_j}{2\gamma}},
\end{equation*}
exists. In particular, the infimum of the functions $\Upsilon(\alpha_{ij},\beta_{ji})$ can be obtained when
\begin{equation}
	(1-\phi_i)\alpha_{ij}=(1-\psi_j)\beta_{ji},
\end{equation}
and
\begin{equation}
	\alpha_{ij}=\beta_{ji}=0, \quad \textrm{as} \quad (1-\phi_i)(1-\psi_j)>\cos_+^2\brac{\frac{x_i-y_j}{2\gamma}}.
\end{equation}
Take
\begin{equation*}
	\phi_i=1-e^{-\phi_i'},\quad
	\psi_j=1-e^{-\psi_j'},
\end{equation*}
in \eqref{eqn:WFR_atomic_con} and drop all the primes, we have
\begin{multline} \label{eqn:WFR_atomic_con2}
	\WFR^2\brac{\sum_{i = 1}^N \mu_i \delta_{x_i},\sum_{j = 1}^M \nu_j \delta_{y_j}}= \\
		2\gamma^2\sup_{\phi_i,\psi_j}\cbra{\sum_i(1-e^{-\phi_i})\mu_i+\sum_j(1-e^{-\psi_j})\nu_j\;:\; \phi_i+\psi_j\le c(x_i,y_j)},
\end{multline}
where
\begin{equation*}
	c(x,y) = -\log\brac{\cos_+^2\brac{\frac{x_i-y_j}{2\gamma}}}.
\end{equation*}
Direct solving \eqref{eqn:WFR_atomic_con2} is a bit complicated, instead we can consider a simpler optimization problem as follows:
\begin{equation} \label{eqn:WFR_atomic_regul}
	2\gamma^2\sup_{\phi_i,\psi_j}\cbra{\sum_i(1-e^{-\phi_i})\mu_i+\sum_j(1-e^{-\psi_j})\nu_j+\epsilon\sum_{ij}\brac{1-e^{\frac{\phi_i+\psi_j}{\epsilon}}}e^{-\frac{c_{ij}}{\epsilon}}},
\end{equation}
with $c_{ij}=c(x_i,y_j)$. This problem is strictly convex by introducing the regularization term. Moreover, the solution of the regularized problem \eqref{eqn:WFR_atomic_regul} is unique and should converge to the original problem \eqref{eqn:WFR_atomic_con2} as $\epsilon\to0$ \cite{ChPeScVi:18b}.

The regularized problem \eqref{eqn:WFR_atomic_regul} can be easily solved by the iterative proportional fitting procedure \cite{BeCaCuNePe:15, ChPeScVi:18b, Sa:15}. The global maximum of this problem can be obtained by alternatively maximizing the objective function w.r.t. to $\phi_i$ and $\psi_j$, i.e.
\begin{equation}
	\left\{\begin{array}{l}
		\phi^{(\ell+1)}= \arg \max_{\phi}\sum_i\brac{(1-e^{-\phi_i})\mu_i-\epsilon e^{\phi_i/\epsilon}\sum_je^{(\psi_j^{(\ell)}-c_{ij})/\epsilon}}, \\
		\psi^{(\ell+1)}= \arg \max_{\psi}\sum_j\brac{(1-e^{-\psi_j})\nu_j-\epsilon e^{\psi_j/\epsilon}\sum_ie^{(\phi_i^{(\ell+1)}-c_{ij})/\epsilon}}.
	\end{array}\right.
\end{equation}
These maximization steps can be done by pointwise computation:
\begin{equation*}
	\left\{\begin{array}{l}
		\widetilde{\phi}_i^{(\ell+1)}=\brac{\mu_i\big/\brac{\sum_je^{-c_{ij}/\epsilon}\widetilde{\psi}_j^{(\ell)}}}^{1/(1+\epsilon)}, \\
		\widetilde{\psi}_j^{(\ell+1)}=\brac{\nu_j\big/\brac{\sum_ie^{-c_{ij}/\epsilon}\widetilde{\phi}_i^{(\ell+1)}}}^{1/(1+\epsilon)},
	\end{array}\right.
\end{equation*}
in which
\begin{equation*}
	\widetilde{\phi}_i^{(\ell)}=e^{\phi_i^{(\ell)}/\epsilon} ,\quad \widetilde{\psi}_j^{(\ell)}=e^{\psi_j^{(\ell)}/\epsilon}.
\end{equation*}
Hence, we have obtained an efficient numerical method for the WFR metric \eqref{eqn:WFR_atomic_con}. The detailed proof of its convergence can be found in \cite{ChPeScVi:18b}, and we will not repeat here. In the numerical implementation, both the distortion introduced by regularization and  the total number of the iterations depend on the choice of $\epsilon$. Therefore it leads to a trade-off between numerical accuracy and convergence speed, which is not unusual in the numerical world. At the beginning of the iterations, the parameter ε can be chosen to be relatively large to accelerate the convergence. Next, we keep reducing the parameter ε in order to ensure the numerical accuracy. In \cite{Sc:16}, there is a detailed discussion on the techniques arising in the implementation of the algorithm.

\vskip2mm\begin{remark}
	In this section, we only discuss the numerical procedure of the simplest discrete case. But the theory holds true for general situation \cite{ChPeScVi:18, ChPeScVi:18b}. For example,  equation \eqref{eqn:WFR_atomic_con2} can be written as
	\begin{equation} \label{eqn:WFR_general}
		\WFR^2(\mu,\nu)=2\gamma^2\sup_{\phi,\,\psi}\cbra{\int \brac{1-e^{-\phi}}\rd \mu+\int \brac{1-e^{-\psi}} \rd \nu\;:\; \phi(x)+\psi(y)\le c(x,y)}.
	\end{equation}
\end{remark}\vskip2mm

%*****************************************************************************************
%************* The application to the earthquake location problems ****************
%*****************************************************************************************
\section{Application to the Earthquake Location problems} \label{sec:ael}
We now turn to the application of the WFR metric to the earthquake location problems. This can be simply done by defining the distance $\mcd$ in \eqref{eqn:mist_fun} using the WFR metric \eqref{eqn:WFR}-\eqref{eqn:WFR_con},
\begin{equation} \label{eqn:def_of_mcd}
	\mcd(f(t),g(t))=\WFR \brac{f(t)^2,g(t)^2}.
\end{equation}
In the above equation, we square the seismic data to ensure the non-negativeness. Equations \eqref{eqn:inv_prob}-\eqref{eqn:ricker},  \eqref{eqn:WFR}-\eqref{eqn:WFR_con} and \eqref{eqn:def_of_mcd} provide the mathematical model to the waveform based earthquake location problems with WFR metric. We illustrate the convexity of the objective function
\begin{equation} \label{eqn:obj_optim_fun}
	\Theta(\bxi,\tau)=\sum_r\chi_r(\bxi,\tau)=\sum_r\mcd^2(d_r(t),s(\boeta_r,t))
		=\sum_r\WFR^2\brac{d_r(t)^2,s(\boeta_r,t)^2},
\end{equation}
with respect to the earthquake location $\bxi$ and origin time $\tau$. \medskip

\begin{figure} 
	\centering
	\includegraphics[width=0.70\textwidth, height=0.18\textheight]{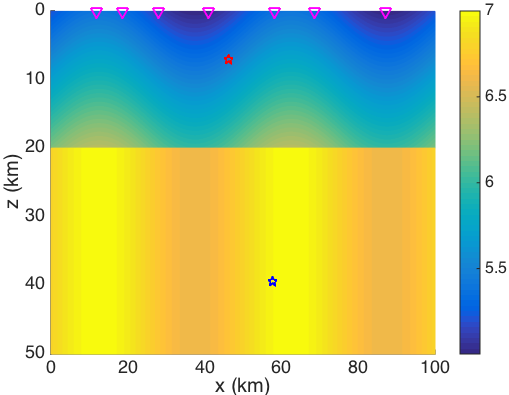}
	\caption{Illustration of two-layer model in Example \ref{exam:two_layler_convexity}. The inverted magenta triangles indicate the receivers. The red and blue pentagram indicate the earthquake hypocenter above and below the Conrad discontinuity respectively.} \label{fig:exam31_vel}
\end{figure}

\begin{figure*} 
	\begin{tabular}{cc}
		\includegraphics[width=0.46\textwidth, height=0.2\textheight]{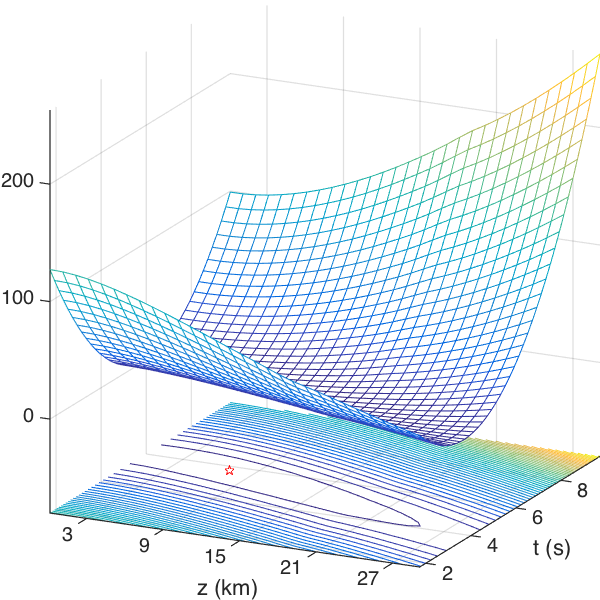} &
		\includegraphics[width=0.46\textwidth, height=0.2\textheight]{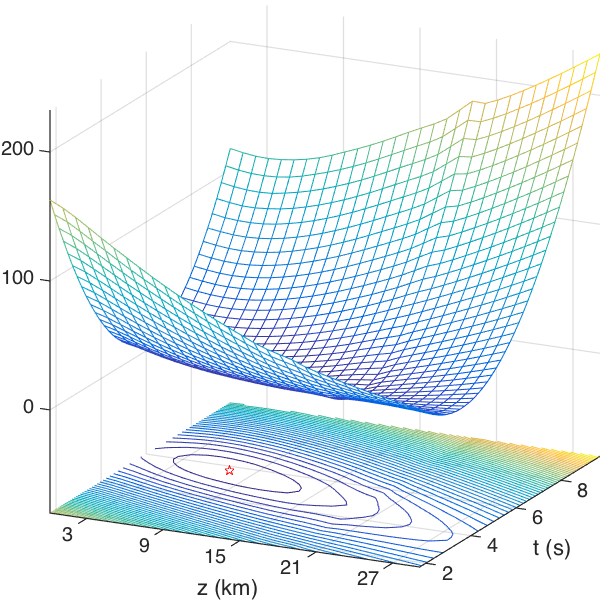} \\
		\includegraphics[width=0.46\textwidth, height=0.2\textheight]{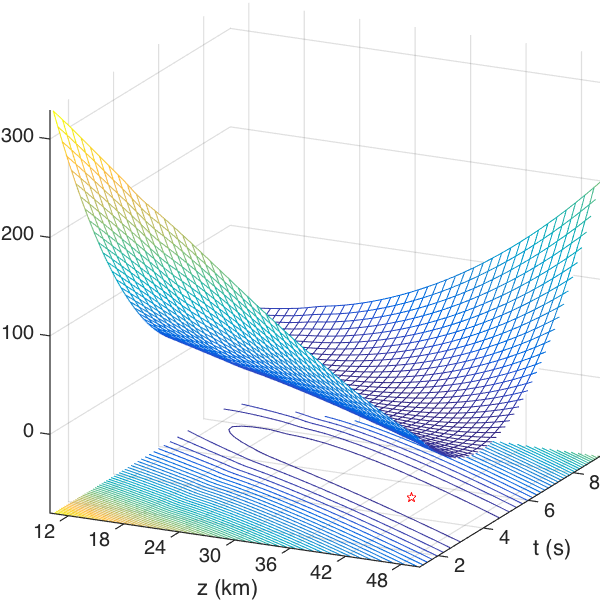} &
		\includegraphics[width=0.46\textwidth, height=0.2\textheight]{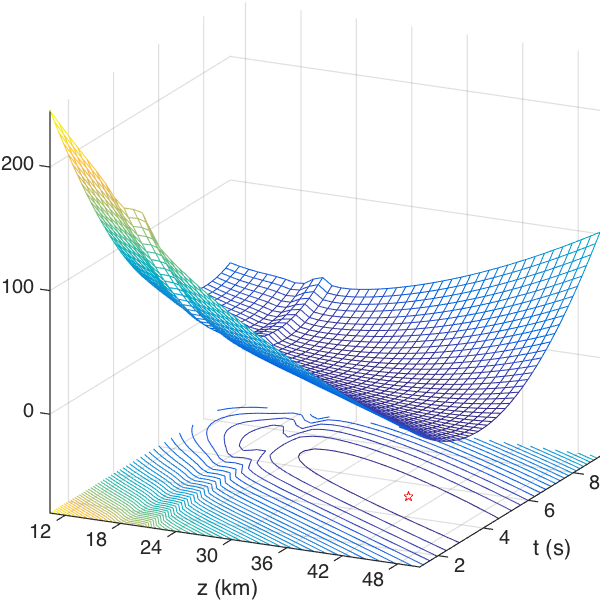}
	\end{tabular}
	\caption{The two-layer model in Example \ref{exam:two_layler_convexity}, the cross-section of the objective function. The red pentagram indicates the minimum point of the objective function, which is also the real earthquake hypocenter. Up for case (i), and Down for case (ii). Left: the quadratic Wasserstein metric; Right: the Wasserstein-Fish-Rao metric.} \label{fig:exam31_IlluCon_NN}
\end{figure*}

\vskip2mm\begin{example} \label{exam:two_layler_convexity}
	This is a 2D two-layer model in the bounded region $\Omega=[0,\,100\,km]\times[0,\,50\,km]$. The Conrad discontinuity is located at a depth of 20km from the Earth's surface, and the wave speed $c(x,z)$ is
	\begin{equation*}
		c(x,z)=\begin{cases}
				5.2+0.05z+0.2sin(\pi x/25), &  z \leq 20\,km, \\
				6.8+0.2sin(\pi x/25), &  z>20\,km. \\
			\end{cases}
	\end{equation*}
	The unit is `km/s'. We randomly set up $7$ receivers on the Earth's surface ($z=0\,km$). Their horizontal positions are listed in Table \ref{tab:two_layler_recepos}, see also Figure \ref{fig:exam31_vel} for the illustration. Here we consider two different earthquake hypocenters, one is above the Conrad discontinuity and the other is below the Conrad discontinuity:
	\begin{align*}
		& (i) \quad \bxi_T^a=(46.23\,km,7.12\,km), \quad \tau_T^a=5.73\,s, \\
		& (ii) \quad \bxi_T^b=(57.60\,km,39.36\,km), \quad \tau_T^b=5.18\,s.
	\end{align*} 
	The dominated frequency of the earthquake is $f_0=2$Hz. We output the cross-section of the objective function defined by the WFR metric \eqref{eqn:obj_optim_fun} and the quadratic Wasserstein metric with the normalized square signals \cite{ChChWuYa:18}. According to the previous discussions, change of  the focal depth would affect the amplitude of the seismogram and the arrival time at a receiver. For the quadratic Wasserstein metric, the amplitude information is ignored by the normalization procedure. Moreover, the arrival time affected by the focal depth could be balanced by adjusting the origin time of the earthquake. Thus, the objective function could be flat on the specific direction of the depth-time plane near the real earthquake hypocenter and origin time, as we can observe in Figure \ref{fig:exam31_IlluCon_NN} (Left). However, the objective function defined using the WFR metric has better convexity near the real earthquake hypocenter and origin time, since the metric takes into account the amplitude information, see Figure \ref{fig:exam31_IlluCon_NN} (Right).
\end{example}

\begin{table*}
	\caption{The two-layer model in Example \ref{exam:two_layler_convexity}: the horizontal positions of receivers, with unit `km'.} \label{tab:two_layler_recepos}
	\begin{center}\begin{tabular}{cccccccc} \hline
		$r$ & $1$ & $2$ & $3$ & $4$ & $5$ & $6$ & $7$ \\ \hline
		$x_r$ & $11.74$ & $18.59$ & $27.88$ & $41.11$ & $58.23$ & $68.50$ & $87.01$ \\ \hline
	\end{tabular}\end{center}
\end{table*}
	
\vskip2mm\begin{remark} \label{rem:bulge}
	In the cross-section of the objective function defined by the WFR metric, we can observe a small bulge near the Conrad discontinuity. It dues to a significant change of seismogram at receiver when the earthquake hypocenter crosses the interface. As we know, the wave is partially transmitted and reflected at the interface. It results in a corresponding change in the seismogram at the receiver. For the quadratic Wasserstein metric, the change is ignored by a normalization procedure, while the WFR metric reflects this change accurately. However, as mentioned in Remark \ref{rem:gamma}, with a larger control parameter $\gamma$, the WFR metric can be more like the classical Wasserstein metric, and the bulge near the Conrad discontinuity can be then eliminated. Away from the interface, we can choose a relatively small parameter $\gamma$ to ensure nice convexity of the WFR metric near the real earthquake hypocenter and origin time.
\end{remark}\vskip2mm

From the above example, we can see the objective function defined by the WFR metric has better convexity near the real earthquake hypocenter and origin time. Thus, we can believe that the minimum point of the objective function defined by the WFR metric has smaller deviation compared to the one defined by the quadratic Wasserstein metric under the same magnitude of noise. This implies higher accuracy of the earthquake location results. Next, we illustrate this perspective with a numerical example. 

\begin{figure*} 
	\begin{tabular}{cc}
		\includegraphics[width=0.46\textwidth, height=0.20\textheight]{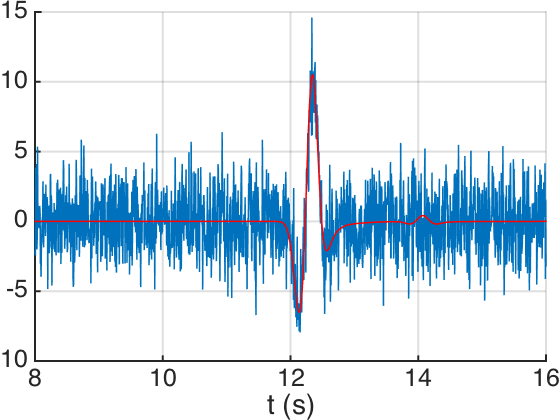} &
		\includegraphics[width=0.46\textwidth, height=0.20\textheight]{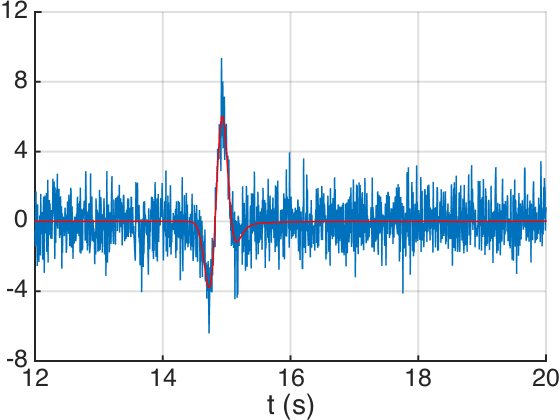}
	\end{tabular}
	\caption{Illustration of signal with noise in the two-layer model. The signal with noise $d_r(t)$ (blue line) and the noise free signal $u(\boeta_r,t;\bxi_T,\tau_T)$ (red line) for receiver $r=1$. The horizontal axis is the time $t$. Up: parameters (i) earthquake hypocenter above the Conrad discontinuity; Down: parameters (ii) earthquake hypocenter below the Conrad discontinuity.} \label{fig:exam32_illu_signal_noise}
\end{figure*}

\begin{figure*} 
	\begin{tabular}{cc}
		\includegraphics[width=0.46\textwidth, height=0.2\textheight]{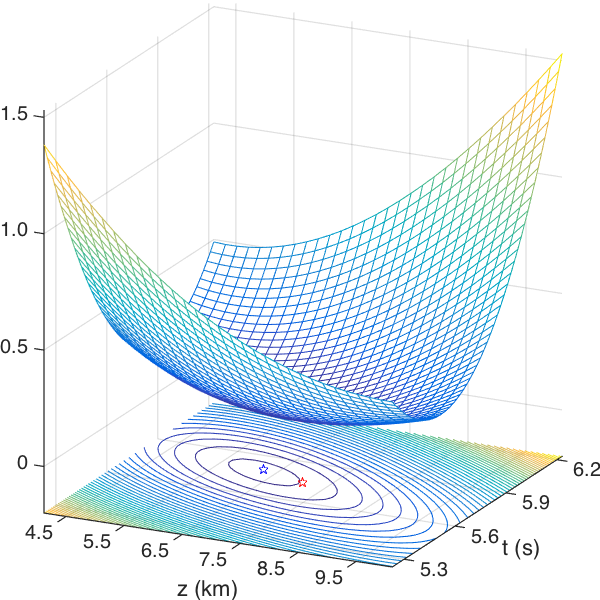} &
		\includegraphics[width=0.46\textwidth, height=0.2\textheight]{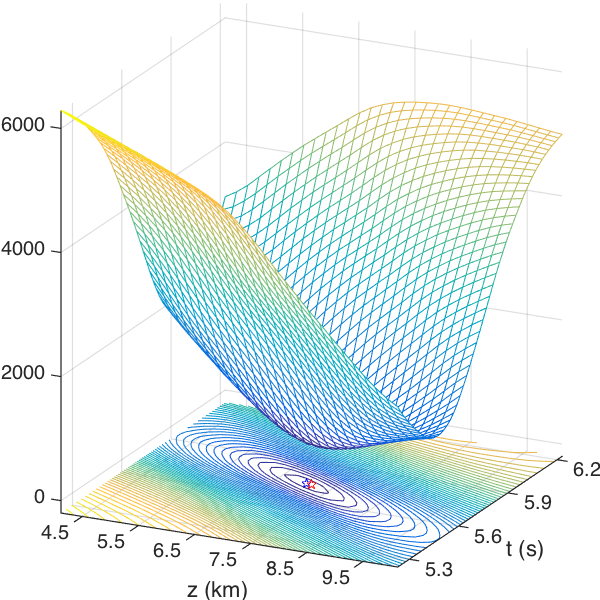} \\
		\includegraphics[width=0.46\textwidth, height=0.2\textheight]{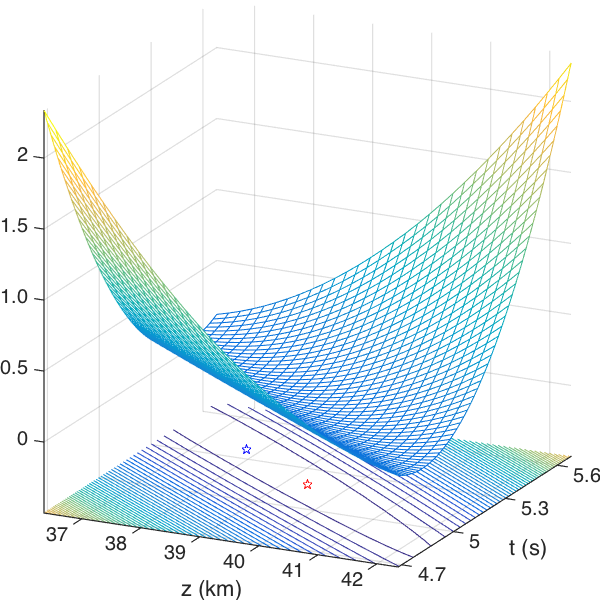} &
		\includegraphics[width=0.46\textwidth, height=0.2\textheight]{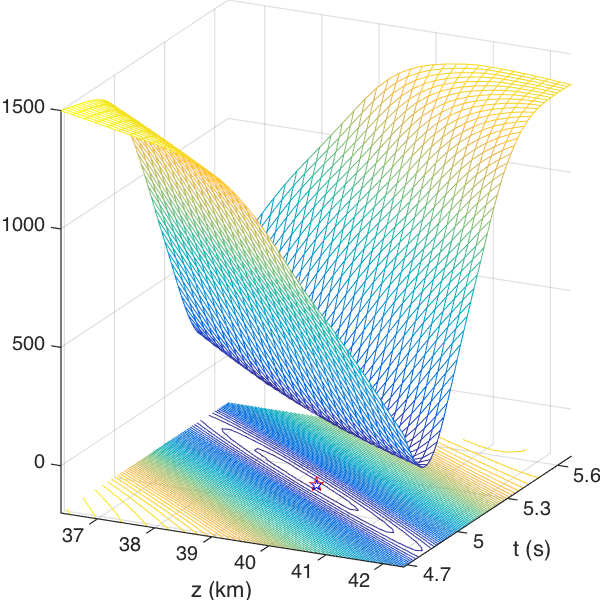}
	\end{tabular}
	\caption{The two-layer model in Example \ref{exam:two_layler_noise_convexity_and_scatter}, the cross-section of the objective function. The blue pentagram denotes the minimum point of the objective function, while the red pentagram denotes the real earthquake hypocenter. Up for case (i), and Down for case (ii). Left: the quadratic Wasserstein metric; Right: the Wasserstein-Fish-Rao metric.} \label{fig:exam32_IlluCon_ND}
\end{figure*}

\begin{figure*} 
	\begin{tabular}{cc}
		\includegraphics[width=0.46\textwidth, height=0.2\textheight]{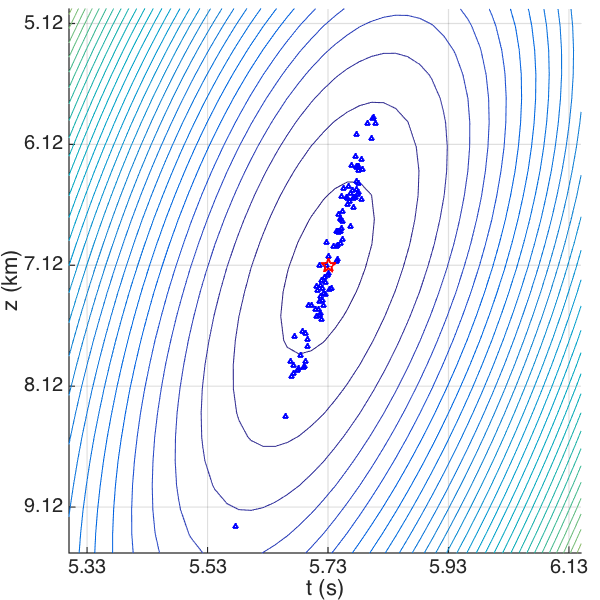} &
		\includegraphics[width=0.46\textwidth, height=0.2\textheight]{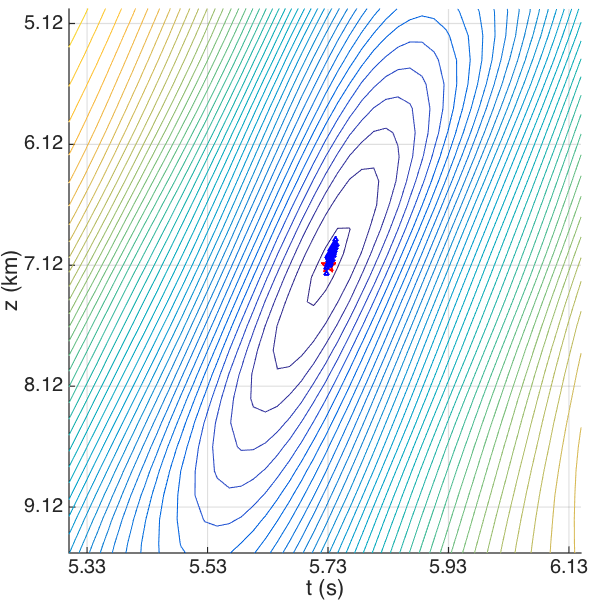} \\
		\includegraphics[width=0.46\textwidth, height=0.2\textheight]{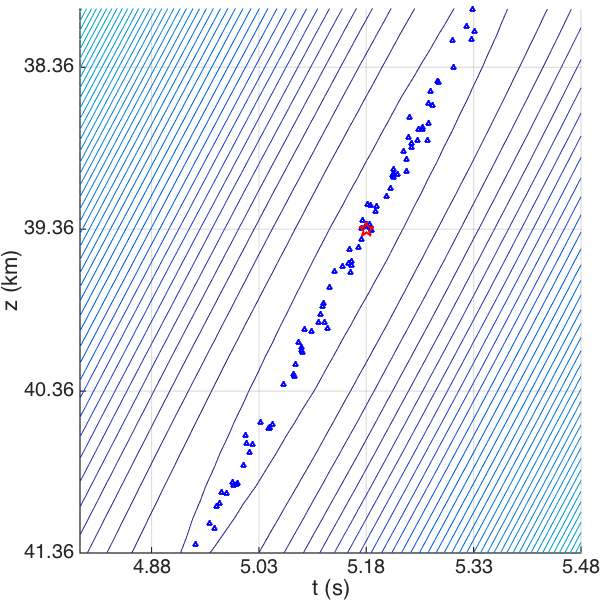} &
		\includegraphics[width=0.46\textwidth, height=0.2\textheight]{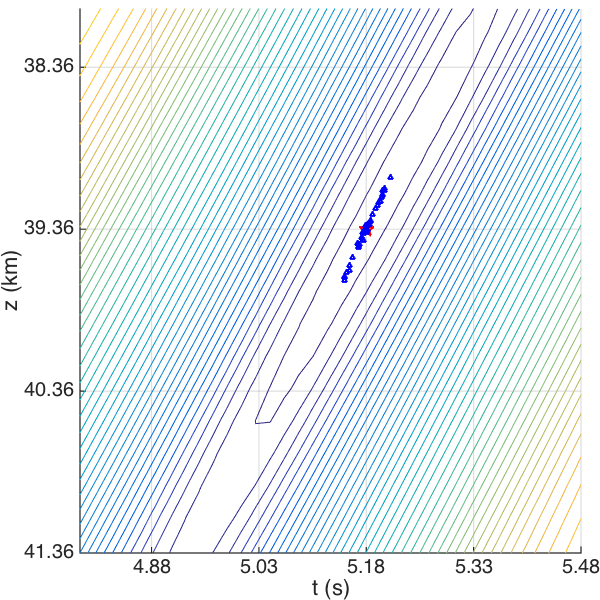}
	\end{tabular}
	\caption{The two-layer model in Example \ref{exam:two_layler_noise_convexity_and_scatter}, the scattergraph of all the minimum points. The blue triangle denotes each minimum points, while the red pentagram denotes the real earthquake hypocenter. The background is a contour plot of the objective function without data noise.Up for case (i), and Down for case (ii). Left: the quadratic Wasserstein metric; Right: the Wasserstein-Fish-Rao metric.} \label{fig:exam32_scatter_ND}
\end{figure*}

\vskip2mm\begin{example} \label{exam:two_layler_noise_convexity_and_scatter}
	Let us consider the same parameters set up as in Example \ref{exam:two_layler_convexity}. Instead of noise-free data, we use the real earthquake signals with noise
	\begin{equation} \label{eqn:add_noise}
		d_r(t)=u(\boeta_r,t;\bxi_T,\tau_T)+N_r(t).
	\end{equation}
	Here $N_r(t)$ is subject to the normal distribution with mean $\mu=0$ and the standard deviation
	\begin{equation} \label{eqn:noise_magnitude}
		\sigma=R\times \max_t\abs{u(\boeta_r,t;\bxi_T,\tau_T)},
	\end{equation}
	in which the ratio $R=20\%$. These signals are illustrated in Figure \ref{fig:exam32_illu_signal_noise}. For the quadratic Wasserstein metric, a time window that contains the main part of $u(\boeta_r,t;\bxi_T,\tau_T)$ is chosen to reduce the impact of noise. For the WFR metric, we do not even need to choose the time windows since the impact of noise is less significant.  This also means that the WFR metric has better robustness when dealing with seismic signals.

	We output the cross-section of the objective function defined by the WFR metric and the quadratic Wasserstein metric with the normalized square signals in Figure \ref{fig:exam32_IlluCon_ND}. As expected, the objective functions in the optimization still have good convexity for both metrics. However, there occur deviations of the minimum point of the objective functions. In Table \ref{tab:two_layler_noise_offset}, the offset between the minimum point and the real earthquake hypocenter are listed for different metrics. from which we can see that the minimum point of the objective function defined by the WFR metric is much closer to the real earthquake hypocenter. Thus, a higher accuracy for the earthquake location can be expected.
	
	Next, we test the effects of noise using $100$ experiments with white noise. We fix the noise ratio $R$ at $20\%$ and choose different random seeds. For each experiment, we simply calculate the all possible values of the objective function through a brute force algorithm and arrive the optimal  point. In each scenario displayed in Figure \ref{fig:exam32_scatter_ND}, all $100$ optimal solution points are presented with background showing the contour of the noise-free objective function. The main purpose of these experiments is to study the propagation of uncertainty in the statistical sense. The averaged and maximum distance between the optimal solution point and the real earthquake hypocenter are given in Table \ref{tab:two_layler_noise_amd}. From the figures and tables, we can conclude that the WFR metric offers a better mathematical tool for the waveform based earthquake location problem.
\end{example} \vskip2mm

\begin{table*}
	\caption{The two-layer model in Example \ref{exam:two_layler_noise_convexity_and_scatter}: the offset between the minimum point and the real earthquake hypocenter for different metrics and different cases, with unit `km'.} \label{tab:two_layler_noise_offset}
	\begin{center}\begin{tabular}{ccc} \hline
		  & quadratic Wasserstein metric & WFR metric \\ \hline
		(i) above the Conrad  & $0.82$ & $0.10$ \\
		(ii) below the Conrad & $1.59$ & $0.02$ \\ \hline
	\end{tabular}\end{center}
\end{table*}

\begin{table*}
	\caption{The two-layer model in Example \ref{exam:two_layler_noise_convexity_and_scatter}: the average(AVE) and maximal(MAX) distance between the minimum point and the real earthquake hypocenter for $100$ tests, with unit `km'. We compare the results of quadratic Wasserstein metric(QWM) and the WFR metric(WFRM).} \label{tab:two_layler_noise_amd}
	\begin{center}\begin{tabular}{ccccc} \hline
		  & AVE-QWM & MAX-QWM & AVE-WFRM & MAX-WFRM \\ \hline
		(i) above the Conrad  & $0.521$ & $2.160$ & $0.072$ & $0.214$ \\
		(ii) below the Conrad & $0.855$ & $2.596$ & $0.130$ & $0.386$ \\ \hline
	\end{tabular}\end{center}
\end{table*}

The above example illustrates that the WFR metric can be used to obtain more accurate earthquake location results than the quadratic Wasserstein metric under the same magnitude of noise. Thus, we believe that the WFR metric could be a better tool for the earthquake location problems.

\vskip2mm\begin{remark}
	For the quadratic Wasserstein metric, the accuracy of the earthquake location results of noisy data can be improved by selecting the time window carefully and increasing the time sampling. However, both techniques are not easily implemented. In contrast, the WFR metric can achieve more accurate earthquake location results without the techniques mentioned above.
\end{remark}\vskip2mm

\subsection{The adjoint method} \label{subsec:adjoint}
Up to now, we have proposed the waveform based earthquake location model with WFR metric. To solve this optimization problem, it is necessary to obtain the gradient of the misfit function $\chi_r(\bxi,\tau)$ in \eqref{eqn:inv_prob} under the WFR metric \eqref{eqn:WFR}-\eqref{eqn:WFR_con} with respect to the earthquake hypocenter $\bxi$ and the origin time $\tau$. First, we define
\begin{align*}
	& F(\phi,\psi,f,g)=2\gamma^2 \brac{\int_0^T \brac{1-e^{-\phi(t)}}f^2(t)\rd t + \int_0^T \brac{1-e^{- \psi(t)}}g^2(t)\rd t}, \\
	& \Pi=\cbra{(\phi,\psi) \,:\, \phi(t_1) + \psi(t_2) \leq c(t_1,t_2), \; \forall t_1,t_2\in[0,T]},
\end{align*}
according to \eqref{eqn:WFR_general} and \eqref{eqn:def_of_mcd}, we have
\begin{equation*}
	\WFR^2(f(t)^2,g(t)^2)=\sup_{(\phi,\psi)\in\Pi}F(\phi,\psi,f,g)=F(\phi^*,\psi^*,f,g),
\end{equation*}
where $\phi^*(t)$ and $\psi^*(t)$ are specific functions that ensure the functional $F(\phi,\psi,f,g)$ to reach the supreme. Through a simple discussion, we can get
\begin{equation*}
	\nabla_g\WFR^2(f^2,g^2)=\partial_gF(\phi^*,\psi^*,f,g),
\end{equation*}
thus the perturbation of the WFR metric with respect to $\delta g$ can be written as
\begin{multline}
	\WFR^2\brac{f^2,(g+\delta g)^2}-\WFR^2\brac{f^2,g^2} \\
		=4\gamma^2\int_0^T\brac{1-e^{-\psi^*(t)}}g(t)\delta g(t)\rd t+O(\delta g(t)^2).
\end{multline}

Next, following the similar manner in Section 3.1 of \cite{ChChWuYa:18}, the relation between the perturbation of misfit function $\delta \chi_r$ and the perturbation of earthquake hypocenter $\delta \bxi$ and the origin time $\delta \tau$ can be obtained as
\begin{eqnarray*}
	\delta \chi_r &=& K_r^{\bxi}\cdot\delta\bxi+K_r^{\tau}\delta\tau, \\
	K_r^{\bxi} &=& \int_0^T R(t-\tau)\nabla w_r(\bxi,t)\rd t, \\
	K_r^{\tau} &=& -\int_0^T R'(t-\tau)w_r(\bxi,t)\rd t,
\end{eqnarray*}
where the adjoint wavefield $w_r(\bxi,t)$ satisfies
\begin{align*} 
	& \frac{\partial^2 w_r(\bx,t)}{\partial t^2}=\nabla\cdot\brac{c^2(\bx)\nabla w_r(\bx,t)}+4\gamma^2 (1 - e^{- \psi^*(t)})s(\boeta_r,t)\delta(\bx-\boeta_r) , \quad \bx\in\Omega, \\
	& w_r(\bx,T)=\partial_t w_r(\bx,T)=0, \quad \bx\in\Omega, \\
	& \bn\cdot \brac{c^2(\bx)\nabla w_r(\bx,t)}=0, \quad \bx\in\partial \Omega.
\end{align*}

% ********************************************************************************
% ************ Numerical Experiments ****************************************
% ********************************************************************************
\section{Numerical Experiments}\label{sec:numer}
In this section, we present two numerical experiments to demonstrate the validity of the WFR metric. We first output the convergence trajectories in the absence of noise. From which, the convergence of the new model defined by WFR metric can be validated. For more important earthquake location problem with data noise, the convergence trajectories are also output. In addition, we will compare the earthquake location results by the quadratic Wasserstein metric and the WFR metric.

In all the numerical examples, the finite difference schemes \cite{Da:86, LiYaWuMa:17, YaLuWuPe:04} are used to solve the acoustic wave equation \eqref{eqn:wave}. The perfectly matched layer boundary condition \cite{KoTr:03} is applied inside the earth, and the reflection boundary condition \eqref{bc:wave} is used on the earth's surface. The point source $\delta(\bx-\bxi)$ is discretized by fifth-order piecewise polynomials \cite{We:08},
\begin{equation*}
	\delta(x)=\left\{\begin{array}{ll}
		\frac{1}{h}\left(1-\frac{5}{4}\abs{\frac{x}{h}}^2-\frac{35}{12}\abs{\frac{x}{h}}^3
			+\frac{21}{4}\abs{\frac{x}{h}}^4-\frac{25}{12}\abs{\frac{x}{h}}^5\right), & \abs{x}\le h, \\
		\frac{1}{h}\left(-4+\frac{75}{4}\abs{\frac{x}{h}}-\frac{245}{8}\abs{\frac{x}{h}}^2+\frac{545}{24}\abs{\frac{x}{h}}^3
			-\frac{63}{8}\abs{\frac{x}{h}}^4+\frac{25}{24}\abs{\frac{x}{h}}^5\right), & h<\abs{x}\le 2h, \\
		\frac{1}{h}\left(18-\frac{153}{4}\abs{\frac{x}{h}}+\frac{255}{8}\abs{\frac{x}{h}}^2-\frac{313}{24}\abs{\frac{x}{h}}^3
			+\frac{21}{8}\abs{\frac{x}{h}}^4-\frac{5}{24}\abs{\frac{x}{h}}^5\right), & 2h<\abs{x}\le 3h, \\
		0, & \abs{x}>3h.
	\end{array}\right. 
\end{equation*}
Here $h$ is a numerical parameter related to the mesh size.

\subsection{The two-layer model} \label{subsec:two_layer}

\begin{figure*} 
	\begin{tabular}{ccc}
		\includegraphics[width=0.3\textwidth, height=0.2\textheight]{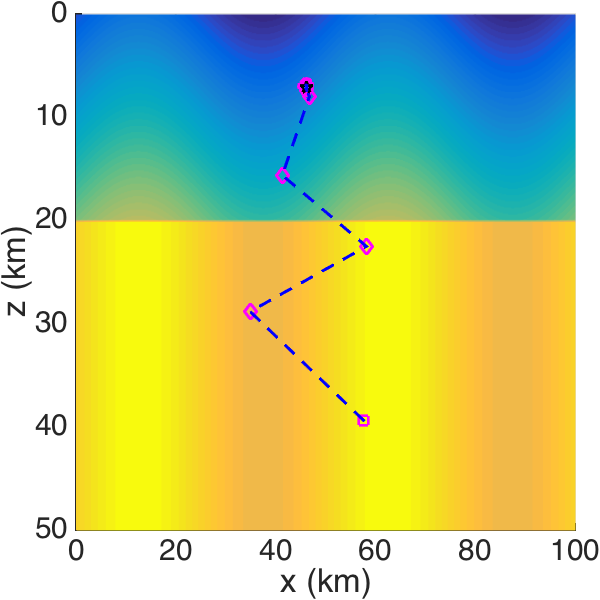} &
		\includegraphics[width=0.3\textwidth, height=0.2\textheight]{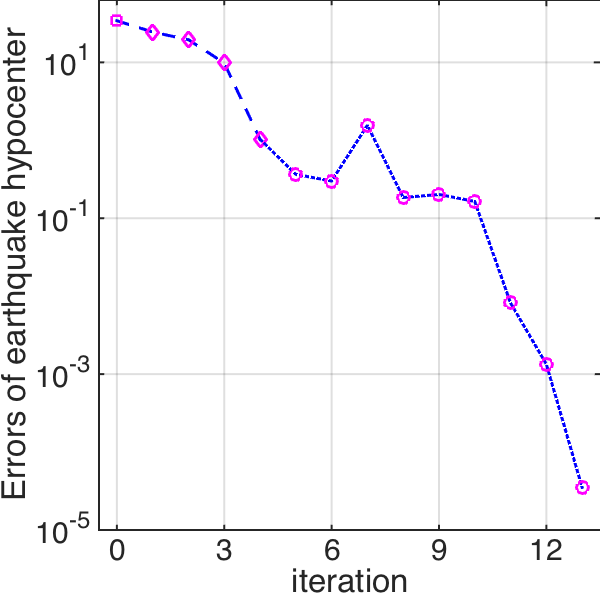} &
		\includegraphics[width=0.3\textwidth, height=0.2\textheight]{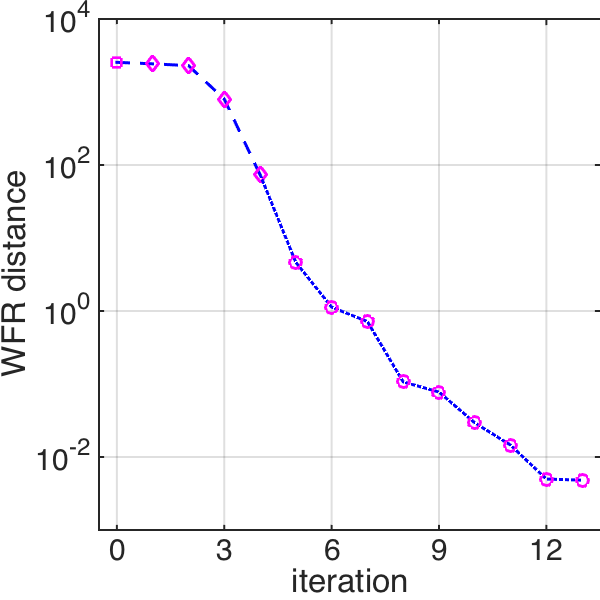} \\
		\includegraphics[width=0.3\textwidth, height=0.2\textheight]{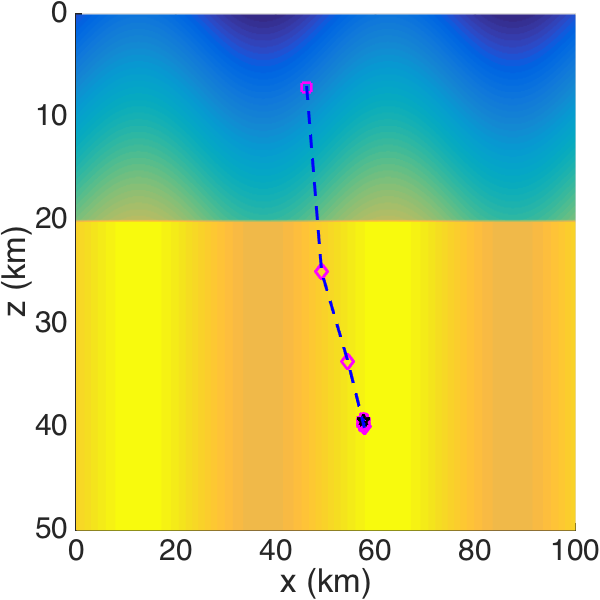} &
		\includegraphics[width=0.3\textwidth, height=0.2\textheight]{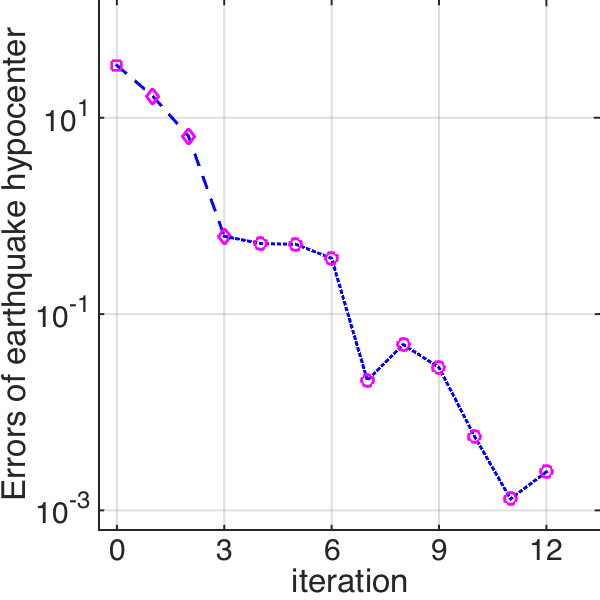} &
		\includegraphics[width=0.3\textwidth, height=0.2\textheight]{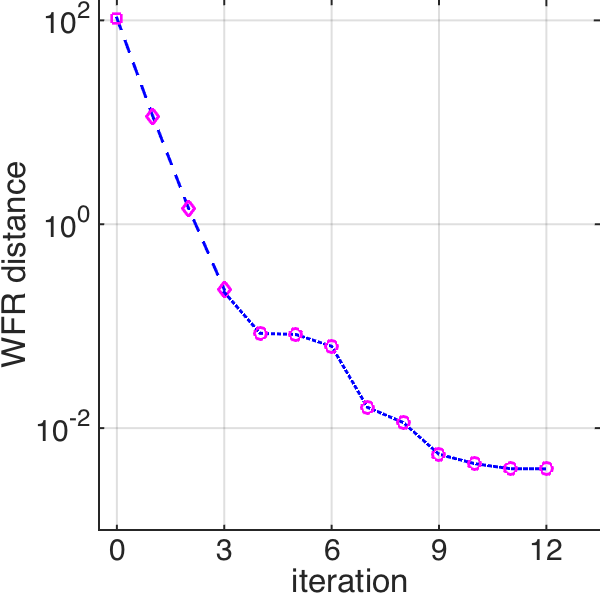}
	\end{tabular}
	\caption{Convergence history of the two-layer model.  Up for case (i), and Down for case (ii). Left: the convergent trajectories; Mid: the absolute errors between the real and computed earthquake hypocenter with respect to iteration steps; Right: the WFR distance between the real and synthetic earthquake signals with respect to iteration steps. The magenta square is the initial hypocenter, the magenta diamond and circle denote the hypocenter in the iterative process with different control parameter $\gamma_i$, and the black pentagram is the real hypocenter.} \label{fig:exam41_traj}
\end{figure*}

Consider the same parameters set up as in Example \ref{exam:two_layler_convexity}, the velocity model is illustrated in Figure \ref{fig:exam31_vel}. We first study the convergency of the new model defined by WFR metric under the ideal situation that there is no noise in the data. We consider the case that the earthquake occurs above the Conrad discontinuity but the initial hypocenter of the earthquake is chosen below the Conrad discontinuity (i), and its contrary case (ii).
\begin{align*}
	&(i) \; \bxi_T=(46.23km, \,7.12km),\;\tau_T=5.73s, \quad \bxi=(57.60km,\,39.36km) ,\tau=3.18s,  \\
	&(ii) \;  \bxi_T=(57.60km,\, 39.36km) ,\;\tau_T=3.18s, \quad \bxi=(46.23km,\, 7.12km), \tau=5.73s.
\end{align*}
In Figure \ref{fig:exam41_traj}, we can see the convergent trajectories, the absolute errors of the earthquake hypocenter and the WFR distance.  As mentioned in Remark \ref{rem:bulge}, the iteration is divided into two phases. In the first phase, a large control parameter is selected to ensure the smoothness of the objective function across the interface. In the second phase, a small control parameter is selected to ensure the local convexity of the objective function near the real earthquake hypocenter and origin time. These figures show that the new method could converge to the real earthquake hypocenter and origin time.

%  \gamma_1=133 or 0.93 ;\gamma_2=0.19

\vskip2mm\begin{remark}
	In Figure \ref{fig:exam41_traj}, we can observe that the WFR distance always decreases, but the error between the real earthquake hypocenter and the hypocenter computed during the iteration  does not keep decreasing. This is because the contour of the WFR distance approximates a very flat ellipse, see Figure \ref{fig:exam31_IlluCon_NN} for illustration. However, as the iteration step increases, the WFR distant and the error decrease on a relatively long scale.
\end{remark}\vskip2mm

\begin{figure*} 
	\begin{tabular}{ccc}
		\includegraphics[width=0.3\textwidth, height=0.2\textheight]{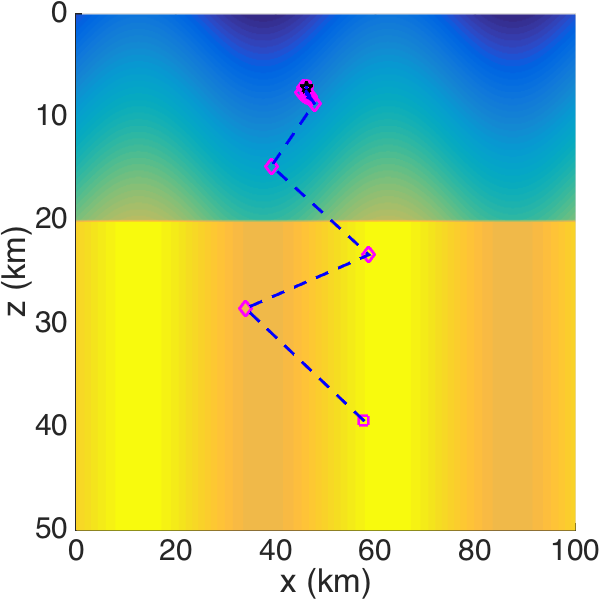} &
		\includegraphics[width=0.3\textwidth, height=0.2\textheight]{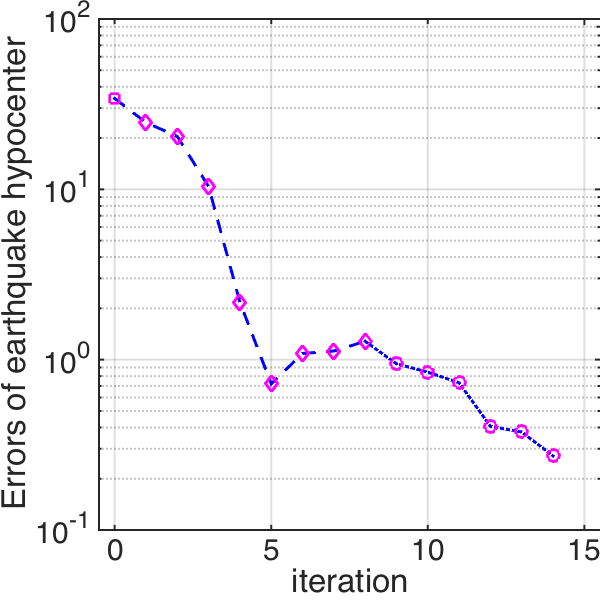} &
		\includegraphics[width=0.3\textwidth, height=0.2\textheight]{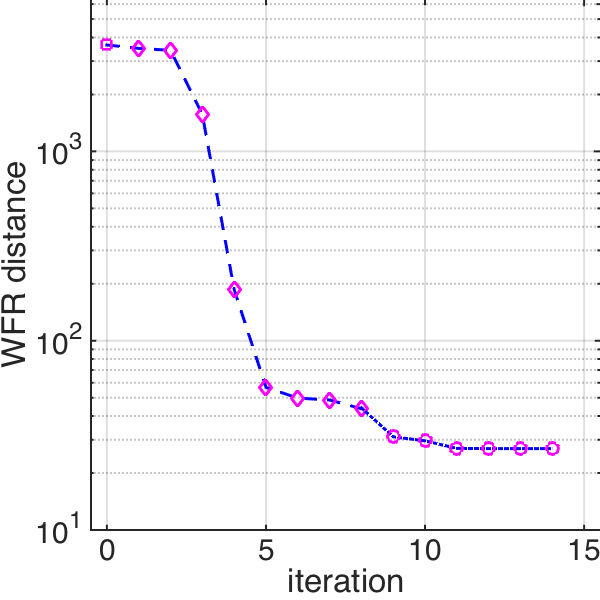} \\
		\includegraphics[width=0.3\textwidth, height=0.2\textheight]{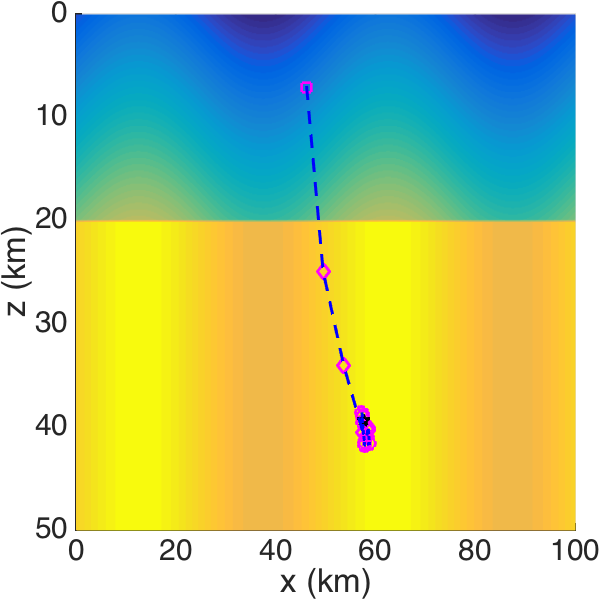} &
		\includegraphics[width=0.3\textwidth, height=0.2\textheight]{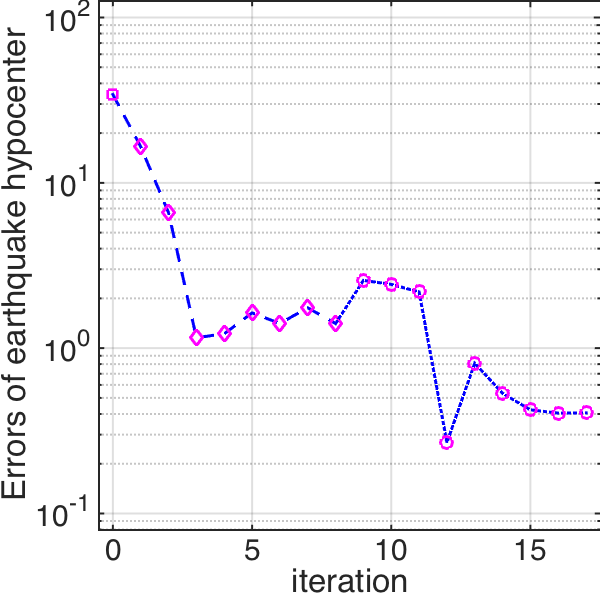} &
		\includegraphics[width=0.3\textwidth, height=0.2\textheight]{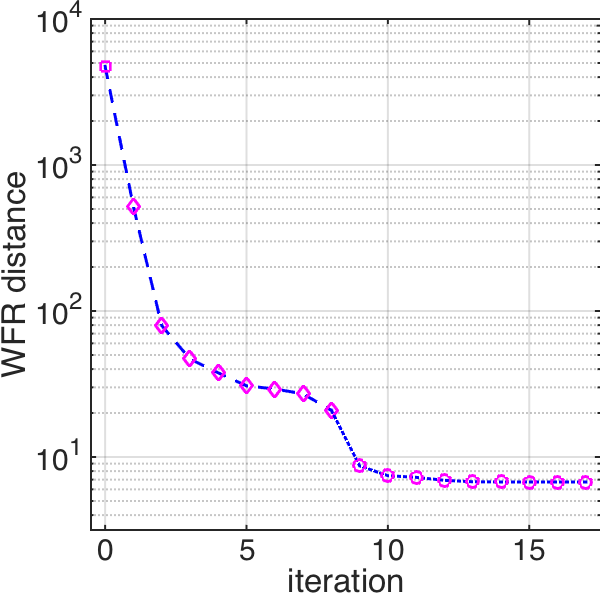}
	\end{tabular}
	\caption{Convergence history of the two-layer model with noise data.  Up for case (i), and Down for case (ii). Left: the convergent trajectories; Mid: the absolute errors between the real and computed earthquake hypocenter with respect to iteration steps; Right: the WFR distance between the real and synthetic earthquake signals with respect to iteration steps. The magenta square is the initial hypocenter, the magenta diamond and circle denote the hypocenter in the iterative process with different control parameter $\gamma_i$, and the black pentagram is the real hypocenter.} \label{fig:exam41_noise_traj}
\end{figure*}

Next, we investigate the adaptability to the noise of the new method. As the equations \eqref{eqn:add_noise}-\eqref{eqn:noise_magnitude} in Example \ref{exam:two_layler_noise_convexity_and_scatter}, noise is added to the real earthquake signal. We select a relatively strong noise $R=20\%$, see Figure \ref{fig:exam32_illu_signal_noise} for illustration. The convergent histories are output in Figure \ref{fig:exam41_noise_traj}. In Table \ref{tab:exam41_noise}, we compare the location errors of different cases and different metrics.  Based on the discussions in Example \ref{exam:two_layler_convexity} and \ref{exam:two_layler_noise_convexity_and_scatter}, we are not surprised that the new earthquake method based on the WFR metric could obtain higher accuracy of the location results when the data is contaminated with noise.

\begin{table*}
	\caption{The two-layer model: the location errors of the iterative solution, with unit `km'. We compare the results of different cases and different metric.} \label{tab:exam41_noise}
	\begin{center}\begin{tabular}{ccc} \hline
		  & quadratic Wasserstein metric & WFR metric \\ \hline
		case (i)  & $1.28$ & $0.27$ \\
		case (ii) & $1.41$ & $0.41$ \\ \hline
	\end{tabular}\end{center}
\end{table*}

\subsection{The subduction plate model} \label{subsec:subduction}
\begin{figure} 
	\centering
	\includegraphics[width=0.70\textwidth, height=0.3\textheight]{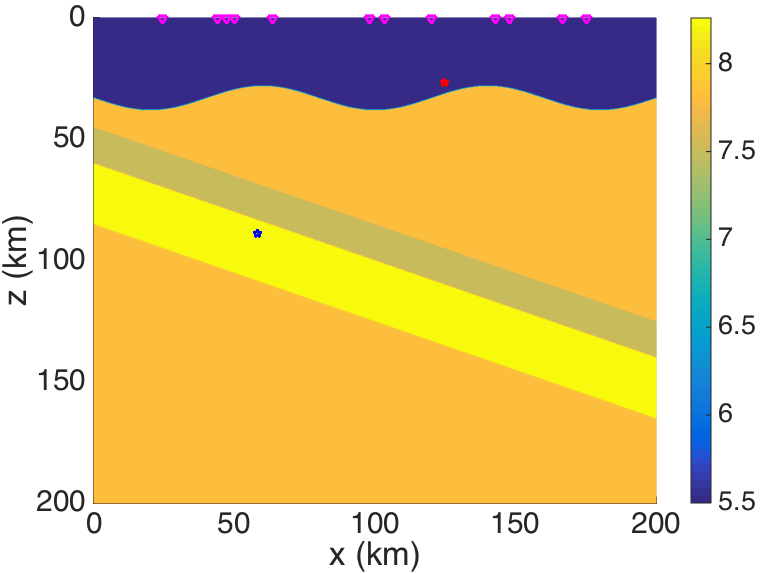}
	\caption{Illustration of subduction plate model. The inverted magenta triangles indicate the receivers. The red and blue pentagram indicate the earthquake hypocenter outside and inside the subduction zone respectively.} \label{fig:exam42_vel}
\end{figure}

Consider a typical seismogenic zone model discussed in \cite{ToZhYa:11}. It consists of the crust, the mantle, and the undulating Moho discontinuity. In addition, there is a subduction zone with a thin low-velocity layer atop a fast velocity layer in the mantle. The earthquake may occur in any of these areas. Taking into account the complex velocity structure, it is much difficult to locate the earthquake. In the simulating domain $\Omega = [0,200\;\mathrm{km}] \times [0,200\;\mathrm{km}]$, the wave speed is
\begin{equation*}
	c(x,z)=\left\{\!\begin{aligned}
		&5.5, & &0 < z \leq 33 + 5 \sin \frac{\pi x}{40}, \\
		&7.8, & &33 + 5 \sin \frac{\pi x}{40} < z \leq 45 + 0.4x, \\
		&7.488, & &45+0.4x < z \leq 60 + 0.4x, \\
		&8.268, & &60+0.4x < z \leq 85 + 0.4x, \\
		&7.8, & &\textrm{others.}
	\end{aligned}\right. 
\end{equation*}
with unit `km/s'. There are 12 randomly distributed receivers $\boldsymbol{\eta}_r = (\chi_r, z_r)$ on the surface $z_r = 0\;\mathrm{km}$. In table \ref{tab:subduction_plate_recepos}, we output their horizontal positions. This velocity model is illustrated in Figure \ref{fig:exam42_vel}. The dominant frequency of the earthquake is $f_0 = 2\; \mathrm{Hz}$ and their simulating time interval $I = [0,55\;\mathrm{s}]$.

\begin{table*}
	\caption{The subduction plate model : the horizontal positions of receivers, with unit `km'.} \label{tab:subduction_plate_recepos}
	\begin{center}\begin{tabular}{ccccccc} \hline
		$r$ & $1$ & $2$ & $3$ & $4$ & $5$ & $6$ \\ \hline
		$x_r$ & $24.53$ & $43.79$ & $47.21$ & $49.96$ & $63.42$ & $97.91$ \\ \hline
		$r$ & $7$ & $8$ & $9$ & $10$ & $11$ & $12$ \\ \hline
		$x_r$ & $103.32$ & $120.11$ & $142.92$ & $147.86$ & $166.73$ & $175.35$ \\ \hline
	\end{tabular}\end{center}
\end{table*}

\begin{figure*} 
	\begin{tabular}{ccc}
		\includegraphics[width=0.3\textwidth, height=0.2\textheight]{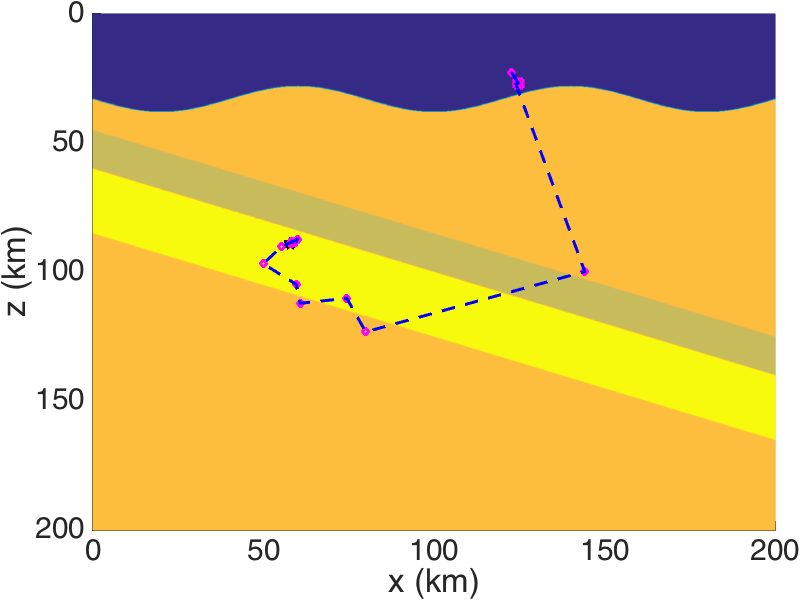} &
		\includegraphics[width=0.3\textwidth, height=0.2\textheight]{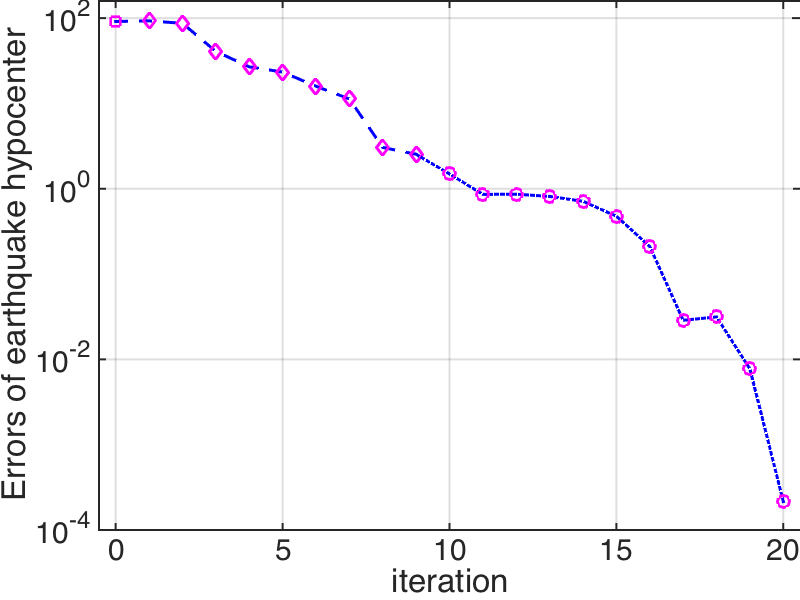} &
		\includegraphics[width=0.3\textwidth, height=0.2\textheight]{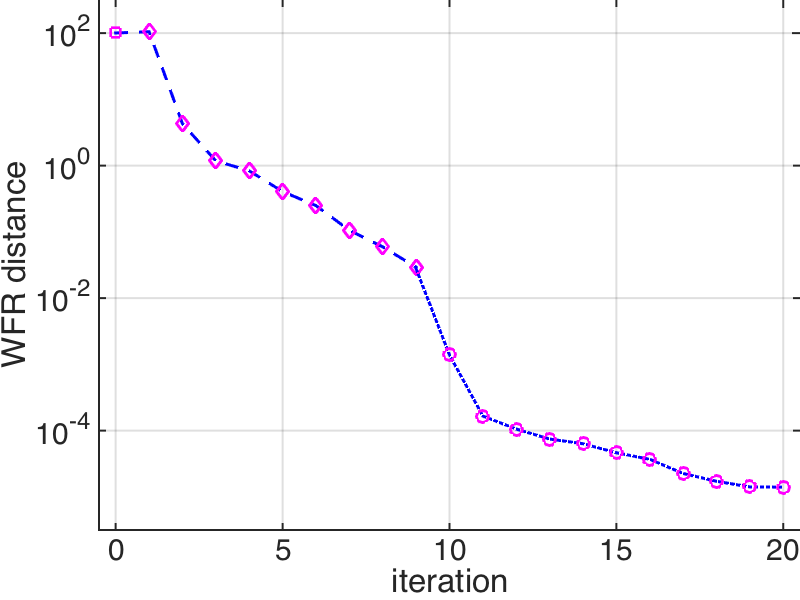} \\
		\includegraphics[width=0.3\textwidth, height=0.2\textheight]{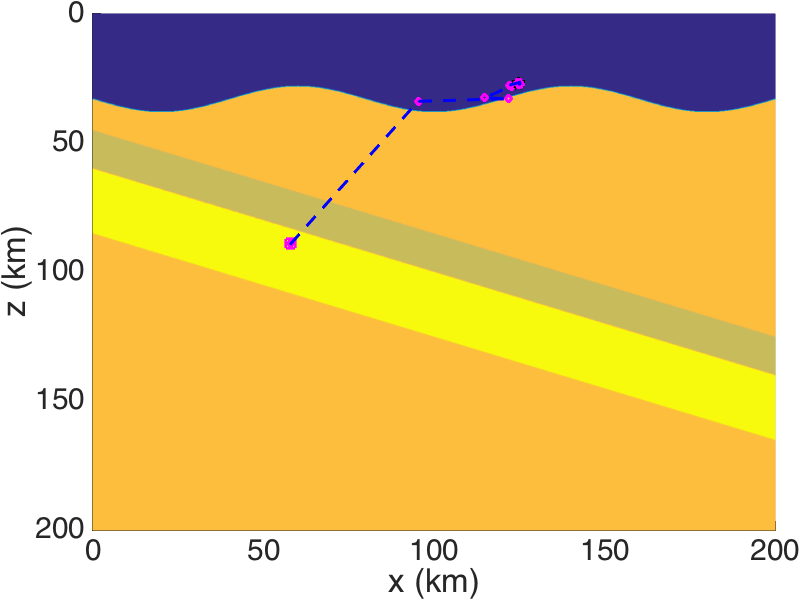} &
		\includegraphics[width=0.3\textwidth, height=0.2\textheight]{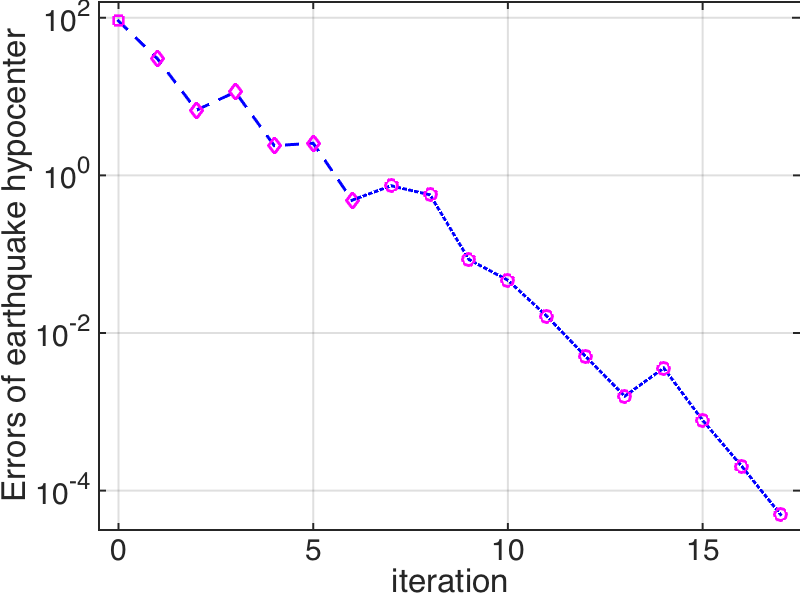} &
		\includegraphics[width=0.3\textwidth, height=0.2\textheight]{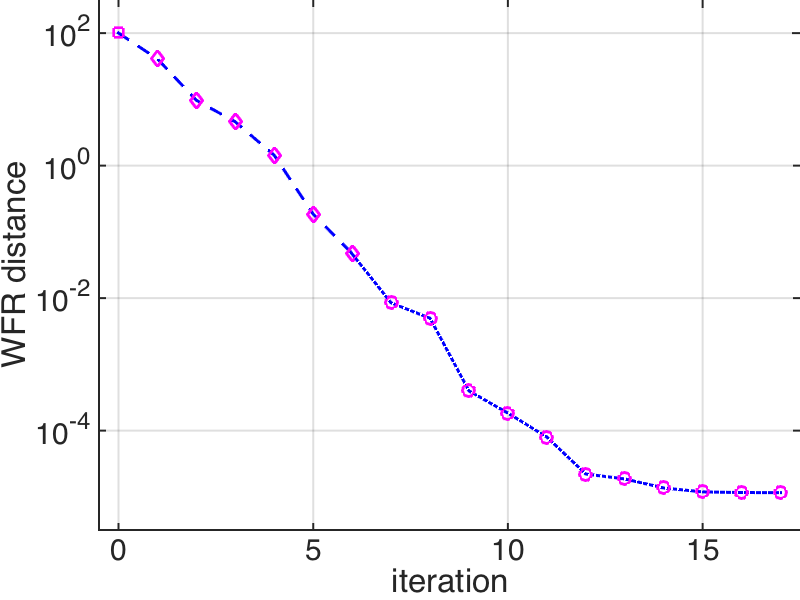}
	\end{tabular}
	\caption{Convergence history of the subduction plate model.  Up for case (i), and Down for case (ii). Left: the convergent trajectories; Mid: the absolute errors between the real and computed earthquake hypocenter with respect to iteration steps; Right: the WFR distance between the real and synthetic earthquake signals with respect to iteration steps. The magenta square is the initial hypocenter, the magenta diamond and circle denote the hypocenter in the iterative process with different control parameter $\gamma_i$, and the black pentagram is the real hypocenter.} \label{fig:exam42_traj}
\end{figure*}

First, consider the ideal situation that there is no noise. We investigate the case when the earthquake occurs in the crust but the initial guess of the earthquake hypocenter is chosen in the subduction zone. Its contrary case is also taken into account. The parameters are selected as follows:
\begin{align*}
	&(i) \; \bxi_T=(124.69km, \,46.48km),\;\tau_T=10.02s, \quad \bxi=(58.96km,\,88.99km) ,\tau=6.79s,  \\
	&(ii) \;  \bxi_T=(58.96km,\,88.99km) ,\;\tau_T=6.79s, \quad \bxi=(124.69km, \,46.48km), \tau=10.02s.
\end{align*}
The convergent trajectories, absolute errors of the earthquake hypocenter and the value of Wasserstein distance are output in Figure \ref{fig:exam42_traj}, from which we can observe nice convergence property of the new method.

\begin{figure*} 
	\begin{tabular}{cc}
		\includegraphics[width=0.46\textwidth, height=0.20\textheight]{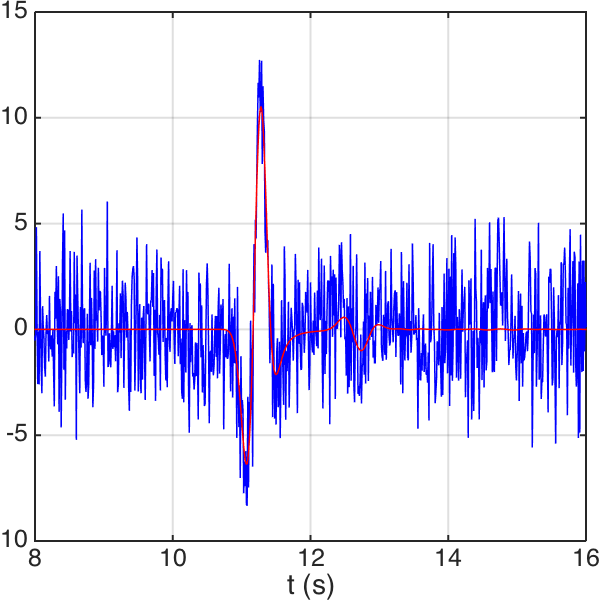} &
		\includegraphics[width=0.46\textwidth, height=0.20\textheight]{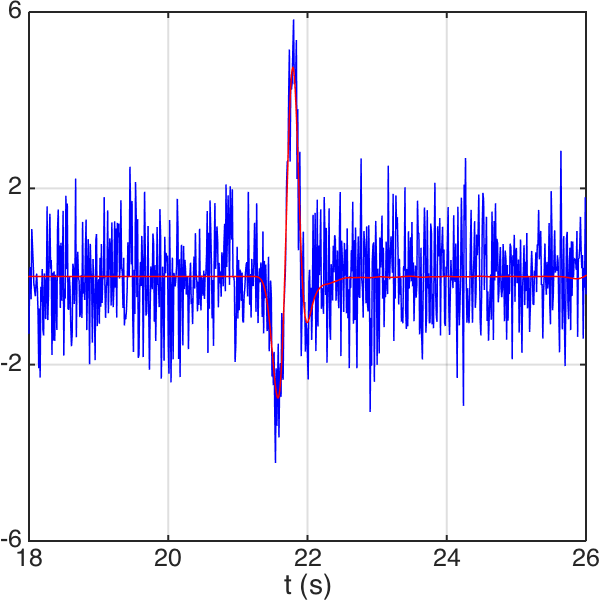}
	\end{tabular}
	\caption{Illustration of signal with noise in the subduction model. The signal with noise $d_r(t)$ (blue line) and the noise free signal $u(\boeta_r,t;\bxi_T,\tau_T)$ (red line) for receiver $r=5$. The horizontal axis is the time $t$. Up: parameters (i) earthquake hypocenter inside the subduction zone; Down: parameters (ii) earthquake hypocenter near the Moho discontinuity.} \label{fig:exam42_illu_signal_noise}
\end{figure*}

\begin{figure*} 
	\begin{tabular}{ccc}
		\includegraphics[width=0.3\textwidth, height=0.2\textheight]{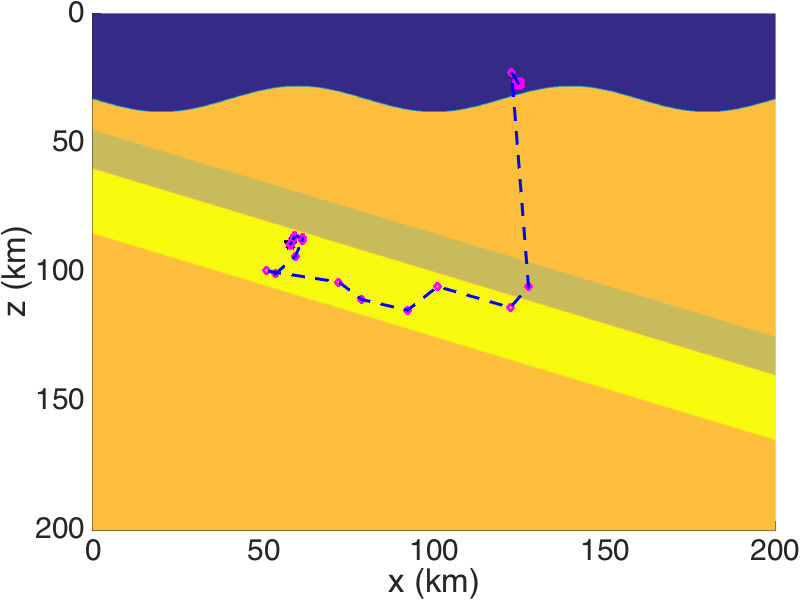} &
		\includegraphics[width=0.3\textwidth, height=0.2\textheight]{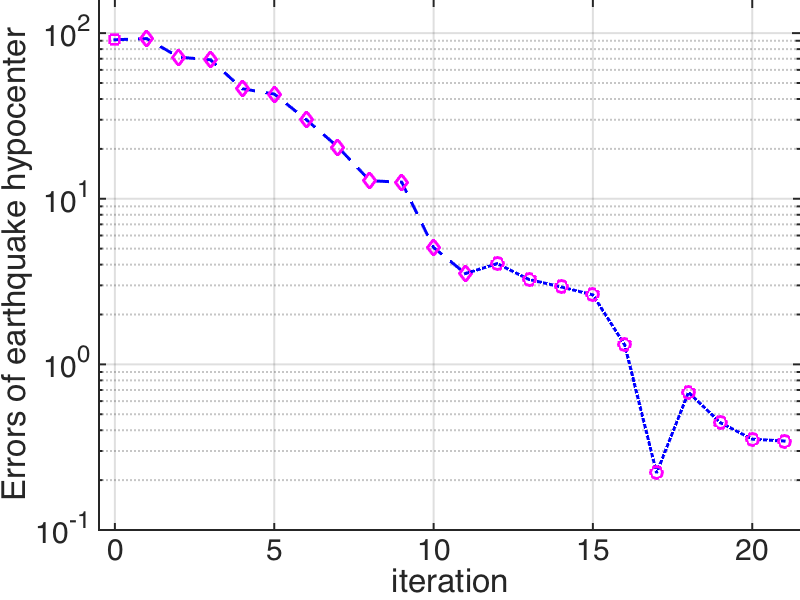} &
		\includegraphics[width=0.3\textwidth, height=0.2\textheight]{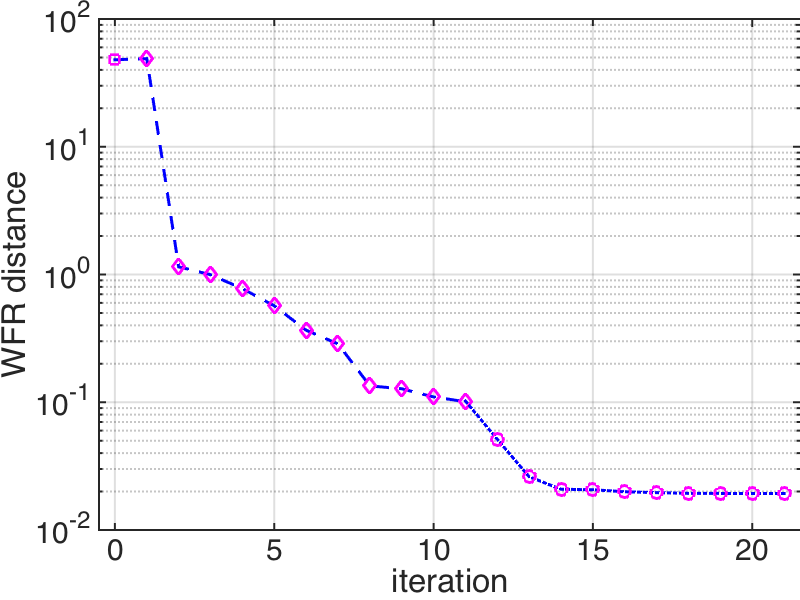} \\
		\includegraphics[width=0.3\textwidth, height=0.2\textheight]{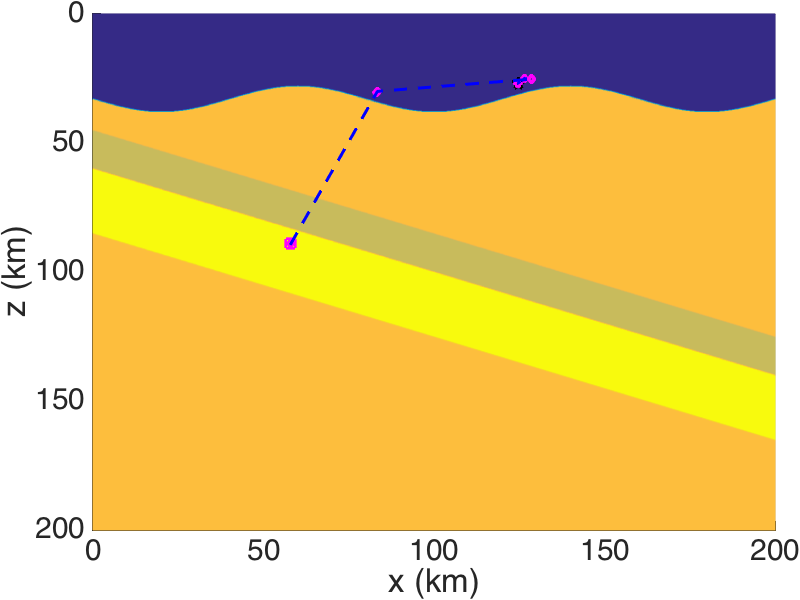} &
		\includegraphics[width=0.3\textwidth, height=0.2\textheight]{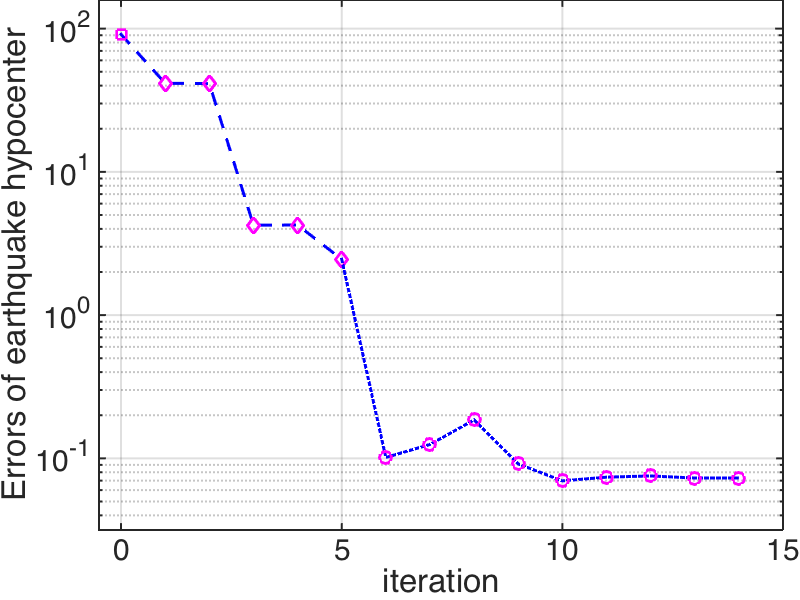} &
		\includegraphics[width=0.3\textwidth, height=0.2\textheight]{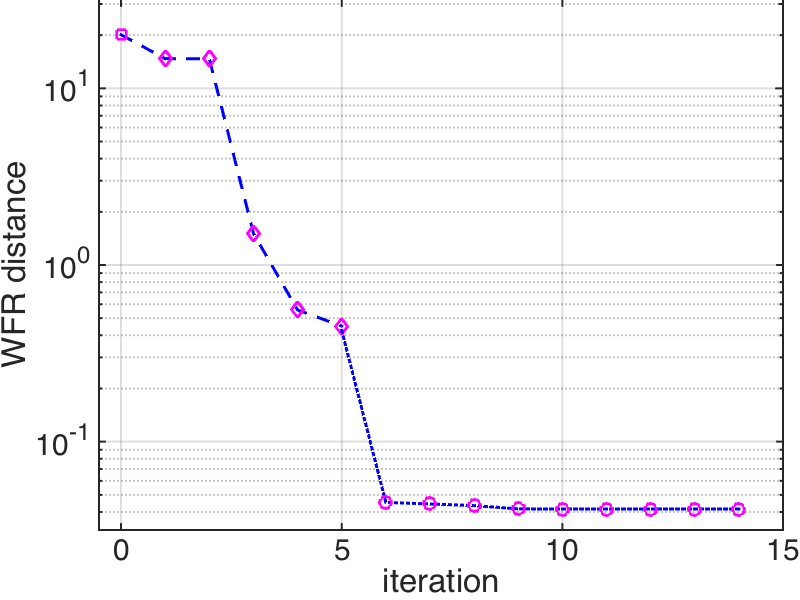}
	\end{tabular}
	\caption{Convergence history of the subduction plate model with noise data.  Up for case (i), and Down for case (ii). Left: the convergent trajectories; Mid: the absolute errors between the real and computed earthquake hypocenter with respect to iteration steps; Right: the WFR distance between the real and synthetic earthquake signals with respect to iteration steps. The magenta square is the initial hypocenter, the magenta diamond and circle denote the hypocenter in the iterative process with different control parameter $\gamma_i$, and the black pentagram is the real hypocenter.} \label{fig:exam42_noise_traj}
\end{figure*}

We also consider the signal containing noise. We select the same parameters (i) and (ii). The noise is added to the real earthquake signals in the same way as in Section \ref{subsec:two_layer}. In Figure \ref{fig:exam42_illu_signal_noise}, we present the real earthquake signal with noise $d_r(t)$ and the noise free signal $u(\boldsymbol{\eta}_r,t;\boldsymbol{\xi}_T,\tau_T)$, and in Figure \ref{fig:exam42_noise_traj} we output the convergent history. In Table \ref{tab:exam42_noise}, we compare the location errors of different cases and different metrics. From which, we can also conclude that the location results of the WFR metric is better than that of the quadratic Wasserstein metric when the data is contaminated with noise.

\begin{table*}
	\caption{The subduction plate model: the location errors of the iterative solution, with unit `km'. We compare the results of different cases and different metric.} \label{tab:exam42_noise}
	\begin{center}\begin{tabular}{ccc} \hline
		  & quadratic Wasserstein metric & WFR metric \\ \hline
		case (i)  & $1.31$ & $0.34$ \\
		case (ii) & $0.54$ & $0.07$ \\ \hline
	\end{tabular}\end{center}
\end{table*}

% ********************************************************************************
% ************ Numerical Experiments *********************************************
% ********************************************************************************
\section{Conclusion}\label{sec:cons}
What we have seen from the above is that the WFR metric is a better choice than the quadratic Wasserstein metric for the earthquake location problems. It overcomes the difficulty that the important amplitude information has been ignored by a normalization procedure in the quadratic Wasserstein metric. By introducing the WFR metric, the retained amplitude information provides more constraints for the earthquake location problems. In particular, the new proposed earthquake location model with WFR metric can achieve more accurate results when the data is contaminated by strong noise.

%*****************************************************************************************
%************* Acknowledgement ******************************************************
%*****************************************************************************************
\section*{Acknowledgement}
The authors are indebted to Prof. Y. Brenier for his helpful discussions.


\begin{thebibliography}{}	
	\bibitem{BeBr:00}
		J.-D. Benamou and Y. Brenier, {\em A computational fluid mechanics solution to the Monge-Kantorovich mass transfer problem}, \textit{Numer. Math.}, \textbf{84}, 375-393, 2000.

	\bibitem{BeCaCuNePe:15}
		J.-D. Benamou, G. Carlier, M. Cuturi, L. Nenna and G. Peyré, {\em Iterative Bregman projections for regularized transportation problems}, \textit{SIAM J. Sci. Comput.}, \textbf{37}(2), A1111-A1138, 2015.		

	\bibitem{ChChWuYa:18}
		J. Chen, Y.F. Chen, H. Wu and D.H. Yang, {\em The quadratic Wasserstein metric for Earthquake Location}, \textit{J. Comput. Phys.}, \textbf{373}, 188-209, 2018.

	\bibitem{ChPeScVi:18}
		L. Chizat, G. Peyr\'e, B. Schmitzer and F.-X. Vialard, {\em An Interpolating Distance Between Optimal Transport and Fisher-Rao Metrics}, \textit{Found. Comput. Math.}, \textbf{18}(1), 1-44, 2018.
		
	\bibitem{ChPeScVi:18b}
		L. Chizat, G. Peyr\'e, B. Schmitzer and F.-X. Vialard, {\em Scaling Algorithms for unbalanced Optimal Transport Problems}, to appear in \textit{Math. Comput.}.

	\bibitem{Da:86}
		M.A. Dablain, The application of high-order differencing to the scalar wave equation, \textit{Geophysics}, \textbf{51}(1), 54-66, 1986.			
		
	\bibitem{EaFo:11}
		D.W. Eaton and F. Forouhideh, {\em Solid angles and the impact of receiver-array geometry on microseismic moment-tensor inversion}, \textit{Geophysics}, \textbf{76}(6): WC77-WC85, 2011.
		
	\bibitem{EnFr:14}
    		B. Engquist and B.D. Froese, {\em Application of the Wasserstein metric to seismic signals}, \textit{Commun. Math. Sci.}, \textbf{12}(5), 979-988, 2014.				

	\bibitem{QiJaAlYaEn:17}
    		L. Qiu, R. Jaime, V. Alejandro, Y. Yang and B. Engquist, {\em Full-waveform inversion with an exponentially encoded optimal-transport norm}, \textit{SEG Technical Program Expanded Abstracts}, 1286-1290, 2017.
		
        \bibitem{EnFrYa:16}
    		B. Engquist, B.D. Froese and Y. Yang, {\em Optimal transport for seismic full waveform inversion}, \textit{Commun. Math. Sci.}, \textbf{14}(8), 2309-2330, 2016.
		
	\bibitem{Fr:12}
		B.D. Froese, \textit{Numerical Methods for the Elliptic Monge-Amp\`ere Equation and Optimal Transport}, Ph.D. thesis, Simon Fraser University, 2012.
		
	\bibitem{Ge:03}
    		M.C. Ge, {\em Analysis of source location algorithms Part I: Overview and non-iterative methods}, \textit{J. Acoust. Emiss.}, \textbf{21}, 14-28. 2003.

	\bibitem{Ge:03b}
    		M.C. Ge, {\em Analysis of source location algorithms Part II: Iterative methods}, \textit{J. Acoust. Emiss.}, \textbf{21}, 29-51, 2003.		
			
	\bibitem{KiLiTr:11}
		Y.H. Kim, Q.Y. Liu and J. Tromp, {\em Adjoint centroid-moment tensor inversions}, \textit{Geophys. J. Int.}, \textbf{186}, 264-278, 2011.
		
	\bibitem{KoMoVo:16}
		S. Kondratyev, L. Monsaingeon and D. Vorotnikov, {\em A new optimal transport distance on the space of finite Radon measures}, \textit{Adv. Differ. Equat.}, \textbf{21}, 1117-1164, 2016.
		
	\bibitem{KoTr:03}
		D. Komatitsch and J. Tromp, {\em A perfectly matched layer absorbing boundary condition for the second-order seismic wave equation}, \textit{Geophys. J. Int.}, \textbf{154}, 146-153, 2003.		
		
	\bibitem{LeSt:81}
		W.H.K. Lee and S.W. Stewart, \textit{Principles and Applications of Microearthquake Networks}, Academic Press, 1981.

	\bibitem{LiYaWuMa:17}
		J.S. Li, D.H. Yang, H. Wu and X. Ma, A low-dispersive method using the high-order stereo-modelling operator for solving 2-D wave equations, \textit{Geophys. J. Int.}, \textbf{210}, 1938-1964, 2017.
	
	\bibitem{LiMiSa:16}
		M. Liero, A. Mielke and G. Savar\'e, {\em Optimal transport in competition with reaction: the Hellinger-Kantorovich distance and geodesic curves}, \textit{SIAM J. Math. Analysis}, \textbf{48}(4), 2869-2911, 2016.
		
	\bibitem{LiMiSa:18}
		M. Liero, A. Mielke and G. Savar\'e, {\em Optimal Entropy-Transport problems and a new Hellinger-Kantorovich distance between positive measures}, \textit{Invent. Math.}, \textbf{211}, 969-1117, 2018.
		
	\bibitem{LiPoKoTr:04}
		Q.Y. Liu, J. Polet, D. Komatitsch and J. Tromp, {\em Spectral-Element Moment Tensor Inversion for Earthquakes in Southern California}, \textit{Bull. seism. Soc. Am.}, \textbf{94}(5), 1748-1761, 2004.
		
	\bibitem{Ma:15}
		R. Madariaga, {\em Seismic Source Theory}, in \textit{Treatise on Geophysics (Second Edition)}, pp. 51-71, ed. Gerald, S., Elsevier B.V., 2015.	
		
	\bibitem{MeBrMeOuVi:16}
		L. M\'etivier, R. Brossier, Q. M\'erigot, E. Oudet and J. Virieux, {\em Measuring the misfit between seismograms using an optimal transport distance: application to full waveform inversion}, \textit{Geophys. J. Int.}, \textbf{205}, 345-377, 2016.
		
	\bibitem{MeBrMeOuVi:16b}
		L. M\'etivier, R. Brossier, Q. M\'erigot, E. Oudet and J. Virieux, {\em An optimal transport approach for seismic tomography: application to 3D full waveform inversion}, \textit{Inverse Probl.}, \textbf{32}, 115008, 2016.
		
	\bibitem{PiRo:14}
		B. Piccoli and F. Rossi, {\em Generalized Wasserstein Distance and its Application to Transport Equations with Source}, \textit{Arch. Rational Mech. Anal. }, \textbf{211}, 335-358, 2014.

	\bibitem{PiRo:16}
		B. Piccoli and F. Rossi, {\em On Properties of the Generalized Wasserstein Distance}, \textit{Arch. Rational Mech. Anal.}, \textbf{222}, 1339-1365, 2016.

	\bibitem{RaPoFi:10}
    		N. Rawlinson, S. Pozgay and S. Fishwick, {\em Seismic tomography: A window into deep Earth}, \textit{Phys. Earth Planet. Inter.}, \textbf{178}, 101-135, 2010.
		
	\bibitem{Ro:67}
		R.T. Rockafellar, {\em Duality and Stability in Extremum Problems Involving Convex Functions}, \textit{Pacific J. Math.}, \textbf{21}(1), 167–187, 1967.

	\bibitem{Sa:15}
		F. Santambrogio, \textit{Optimal Transport for Applied Mathematicians}, Birkh\"auser, 2015.

	\bibitem{SaLoZo:08}
		C. Satriano, A. Lomax and A. Zollo, {\em Real-Time Evolutionary Earthquake Location for Seismic Early Warning}, \textit{Bull. seism. Soc. Am.}, \textbf{98}(3), 1482-1494, 2008.
		
	\bibitem{Sc:16}
		B. Schmitzer, {\em Stabilized Sparse Scaling Algorithms for Entropy Regularized Transport Problems}, \textit{arXiv:1610.06519v}, 2016.
		
	\bibitem{Si:64}
		R. Sinkhorn, {\em A relationship between arbitrary positive matrices and doubly stochastic matrices}, \textit{Ann. Math. Statist.}, \textbf{35}, 876–879, 1964.    

	\bibitem{ToZhYa:11}
    		P. Tong, D.P. Zhao, D.H. Yang, {\em Tomography of the 1995 Kobe earthquake area: comparison of finite-frequency and ray approaches}, \textit{Geophys. J. Int.}, \textbf{187}, 278-302, 2011.

	\bibitem{ToZhYaYaChLi:14}
		P. Tong, D. Zhao, D.H. Yang, X. Yang, J. Chen and Q. Liu, {\em Wave-equation-based travel-time seismic tomography - Part 1: Method}, \textit{Solid Earth}, \textbf{5}, 1151-1168, 2014.

	\bibitem{Vi:03}
		C. Villani, \textit{Topics in Optimal Transportation}, Graduate Studies in Mathematics, American Mathematical Society, 2003.

	\bibitem{Vi:08}
		C. Villani, \textit{Optimal Transport: Old and New}, Springer Science \& Business Media, 2008.

	\bibitem{WaEl:00}
		F. Waldhauser and W.L. Ellsworth, {\em A double-difference earthquake location algorithm: Method and application to the northern Hayward Fault}, California, \textit{Bull. seism. Soc. Am.}, \textbf{90}(6), 1353-1368, 2000.
		
	\bibitem{Wa:04}
		X.-J. Wang, {\em On the design of a reflector antenna II}, \textit{Calc. Var. Partial Dif.}, \textbf{20}(3), 329-341, 2004.
		
	\bibitem{We:08}
    		X. Wen, High Order Numerical Quadratures to One Dimensional Delta Function Integrals, \textit{SIAM J. Sci. Comput.}, \textbf{30}(4), 1825-1846, 2008.				

	\bibitem{WuChHuYa:16}
    		H. Wu, J. Chen, X.Y. Huang and D.H. Yang, {\em A new earthquake location method based on the waveform inversion}, \textit{Commun. Comput. Phys.}, \textbf{23}(1), 118-141, 2018.
					
	\bibitem{WuChJiToYa:17}
		H. Wu, J. Chen, H. Jing, P. Tong and D.H. Yang, {\em The auxiliary function method for waveform based earthquake location}, \textit{arXiv:1706.05551}, 2017.

	\bibitem{WuYa:13}
    		H. Wu and X. Yang, Eulerian Gaussian beam method for high frequency wave propagation in the reduced momentum space, \textit{Wave Motion}, \textbf{50}(6), 1036-1049, 2013.

	\bibitem{YaLuWuPe:04}
		D.H. Yang, M. Lu, R.S. Wu and J.M. Peng, An Optimal Nearly Analytic Discrete Method for 2D Acoustic and Elastic Wave Equations, \textit{Bull. seism. Soc. Am.}, \textbf{94}(5), 1982-1991, 2004.

	\bibitem{YaEnSuFr:18}
    		Y.N. Yang, B. Engquist, J.Z. Sun and B.D. Froese, {\em Application of Optimal transport and the quadratic Wasserstein metric to Full-Waveform-Inversion}, \textit{Geophysics}, \textbf{83}(1), R43-R62, 2018.

	\bibitem{YaEn:18}
    		Y.N. Yang and B. Engquist, {\em Analysis of optimal transport and related misfit functions in full-waveform inversion}, \textit{Geophysics}, \textbf{83}(1), A7-A12, 2018.
		
\end{thebibliography}
\end{document}